\documentclass[12pt]{article}
\usepackage{amsfonts,amssymb,latexsym}

\newtheorem{theorem}{Theorem}[section]
\newtheorem{lemma}[theorem]{Lemma}
\newtheorem{corollary}[theorem]{Corollary}
\newtheorem{proposition}[theorem]{Proposition}
\newenvironment{proof}
{\par\addvspace{0.3cm}\noindent{\rm Proof. }}
{\nopagebreak\mbox{}\hfill $\Box$\par\addvspace{0.25cm}}

\newcommand{\be}{\begin{equation}}
\newcommand{\ee}{\end{equation}}
\newcommand{\bqn}{\begin{eqnarray}}
\newcommand{\eqn}{\end{eqnarray}}
\newcommand{\nn}{\nonumber}
\newcommand{\ba}{\begin{array}}
\newcommand{\ea}{\end{array}}

\newcommand{\under}[2]{\settowidth{\laenge}{$#1$}
\parbox[t]{\laenge}{\makebox[\laenge]{$#1$}\\[-.5ex]
\makebox[\laenge]{$#2$}}}

\newcommand{\cT}{{\cal T}}
\newcommand{\cH}{{\cal H}}
\newcommand{\cM}{{\cal M}}
\newcommand{\cL}{{\cal L}}
\newcommand{\cB}{{\cal B}}
\newcommand{\cN}{{\cal N}}

\newcommand{\C}{{\mathbb C}}
\newcommand{\Z}{{\mathbb Z}}
\newcommand{\T}{{\mathbb T}}
\newcommand{\ka}{\varkappa}
\newcommand{\ro}{\varrho}
\renewcommand{\phi}{\varphi}
\newcommand{\kr}{^\dagger}

\newcommand{\wind}{{\rm wind\,}}
\newcommand{\diag}{{\rm diag\,}}
\newcommand{\ind}{{\rm ind\,}}
\newcommand{\ti}{\tilde}
\newcommand{\wt}{\widetilde}
\newcommand{\wh}{\widehat}
\newcommand{\iv}{^{-1}}
\newcommand{\iy}{\infty}

\newcommand{\NN}{^{N\times N}}
\newcommand{\NNN}{^{2N\times 2N}}
\newcommand{\ovl}{\overline}
\newcommand{\bka}{\ovl{\ka}}

\begin{document}
\newlength{\laenge}

\date{}
\title{Invertibility of Toeplitz + Hankel Operators and Singular Integral
Operators with Flip. -- The case of smooth generating functions}
\author{Torsten Ehrhardt\thanks{tehrhard@mathematik.tu-chemnitz.de.}\\
        	Fakult\"{a}t f\"{u}r Mathematik\\
         	Technische Universit\"{a}t Chemnitz\\
         	09107 Chemnitz, Germany}
\maketitle

\begin{abstract}
It is well known that a Toeplitz operator is invertible if and only if
its symbols admits a canonical Wiener--Hopf factorization, where the
factors satisfy certain conditions.  A similar result holds also for
singular integral operators.  More general, the dimension of the
kernel and cokernel of Toeplitz or singular integral operators which
are Fredholm operators can be expressed in terms of the partial
indices $\ka_1,\dots,\ka_N\in\Z$ of an associated Wiener--Hopf
factorization problem.

In this paper we establish corresponding results for Toeplitz + Hankel
operators and singular integral operators with flip under the
assumption that the generating functions are sufficiently smooth
(e.g., H\"older continuous).  We are led to a slightly different
factorization problem, in which pairs
$(\ro_1,\ka_1),\dots,(\ro_N,\ka_N)\in\{-1,1\}\times\Z$, instead of the
partial indices appear.  These pairs provide the relevant information
about the dimension of the kernel and cokernel and thus answer
the invertibility problem.
\end{abstract}

\section{Introduction}

Let $L^\iy(\T)$ stand for the C*-algebra of all essentially bounded
and Lebesgue measurable functions defined on the unit circle
$\T=\{\;z\in\C\;:\;|z|=1\;\}$, and let $L^2$ stand for the Hilbert
space of all square integrable functions defined on $\T$.  Let $H^2$
($\ovl{H^2}$, resp.) stand for the Hardy space consisting of all
functions $f\in L^2$ for which the Fourier coefficients
\bqn
f_n &=& \frac{1}{2\pi}\int_{0}^{2\pi}f(e^{i\theta})e^{-in\theta}\,d\theta
\eqn
vanish for all $n<0$ ($n>0$, resp.).  Moreover, let
$H^\iy=L^\iy(\T)\cap H^2$ and $\ovl{H^\iy}=L^\iy(\T)\cap\ovl{H^2}$ be
the usual Hardy spaces of essentially bounded function.  Note that
$H^\iy$ and $\ovl{H^\iy}$ are Banach subalgebras of $L^\iy(\T)$.
Finally, let $C(\T)$ stand for the C*-algebra of all continuous
functions defined on the unit circle.

Given a Banach space $X$, let $X^N$ stand for the Banach space of all
$N\times1$ vectors with entries in $X$, and let $X\NN$ stand for the
Banach space of all $N\times N$ matrices with entries in $X$.  Given a
Banach algebra $B$, we denote by $GB$ the group of all invertible
elements in $B$.  A Banach subalgebra $B_1$ of a Banach algebra $B_0$
is called inverse closed (in $B_0$) if $b\in B_1\cap GB_0$ implies
$b\in GB_1$.

For an $N\times N$ matrix valued function $A\in L^\iy(\T)\NN$, the
multiplication operator generated by $A$ is defined by
\bqn
M(A) &:& (L^2)^N\to(L^2)^N,\;\;f(e^{i\theta})\mapsto A(e^{i\theta})
f(e^{i\theta}).
\eqn
The Riesz projection $P$ and the associated projection $Q$
acting on $(L^2)^N$ are given by
\bqn
P\;:\;\sum_{n=-\iy}^\iy f_ne^{in\theta}\mapsto
\sum_{n=0}^\iy f_ne^{in\theta},
&\quad&
Q\;:\;\sum_{n=-\iy}^\iy f_ne^{in\theta}\mapsto
\sum_{n=-\iy}^{-1} f_ne^{in\theta}.
\eqn
Note that $(H^2)^N$ is the image of the Riesz projection $P$.
The flip operator $J$ is defined by
\bqn
J &:& (L^2)^N\to(L^2)^N,\;\;f(e^{i\theta})
\mapsto e^{-i\theta}f(e^{-i\theta}).
\eqn

It is well known that for $A\in L^\iy(\T)\NN$ the Toeplitz operator
\bqn
T(A) &=& PM(A)P
\eqn
acting on $(H^2)^N$ is a Fredholm operator if and only if
$A$ possesses a factorization of the form
\bqn\label{f1.fact}
A(t) &=& A_-(t)\Lambda(t)A_+(t),\qquad t\in\T,
\eqn
where $\Lambda(t)=\diag(t^{\ka_1},\dots,t^{\ka_N})$ is a
diagonal matrix with $\ka_1,\dots,\ka_N\in\Z$, and the factors
$A_+$ and $A_-$ satisfy the following conditions:
\begin{itemize}
\item[(i)]
$A_+\in (H^2)\NN$, $A_+\iv\in (H^2)\NN$;
\item[(ii)]
$A_-\in (\ovl{H^2})\NN$, $A_-\iv\in (\ovl{H^2})\NN$;
\item[(iii)]
The operator $M(A_+\iv)PM(A_-\iv)$, which is a well defined mapping
from $C(\T)^N$ into the Lebesgue space $L^1(\T)^N$, can be extented
by continuity to a linear bounded operator acting from $(L^2)^N$ into
$(L^2)^N$.
\end{itemize}
The integers $\ka_1,\dots,\ka_N$ are called the partial indices
of the above factorization and are uniquely determined up to change of order.
A necessary (but not sufficient) condition for the Fredholmness of
$T(A)$ is that $A\in GL^\iy(\T)\NN$.
If $T(A)$ is a Fredholm operator, then the dimension of
the kernel and cokernel are given by
\begin{equation}
\dim\ker T(A)=-\sum_{\ka_j<0}\ka_j,\qquad
\dim\ker T(A)^*=\sum_{\ka_j>0}\ka_j.
\end{equation}
Here ``$*$'' stands for the adjoint of an operator.
The index of $T(A)$, i.e., the number
$\ind T(A):=\dim\ker T(A)-\dim\ker T(A)^*$, is equal to $-\ka$, where
\bqn
\ka &=&\sum_{j=1}^N\ka_j
\eqn
is the so-called total index of the factorization. In particular, the Toeplitz
operator $T(A)$ is invertible if and only if $A$ admits a canonical
factorization, i.e., a factorization where all partial indices are zero.

A factorization of a matrix function in the form (\ref{f1.fact}) with
the properties (i)--(iii) is sometimes called a generalized
factorization or a $\Phi$-factorization in the space $L^2$.  For
further information about this type of factorization and
generalizations of it, we refer the reader to the monographs
\cite{CG,LS,GK}.

For some classes of functions (e.g., piecewise continuous matrix
functions) there exist different Fredholm criteria, which are easier
to verify.  There also exist explicit formulas for total index.
However, in the case $N>1$, the explicit construction of a
factorization, or, at least the determination of the partial indices
is often the only possibility to answer the question about the
invertibility (and, more general, to calculate the dimension of the
kernel and cokernel in the case of Fredholm operators).

For singular integral operators (with $A,B\in L^\iy(\T)\NN$)
\bqn
S(A,B) &=& PM(A) +QM(B),
\eqn
which are defined on $(L^2)^N$, a similar result holds.  Namely,
$S(A,B)$ is Fredholm if and only if $A,B\in G(L^\iy(\T)\NN)$ and if
$T(AB\iv)$ is a Fredholm operator.  The latter means that the matrix
function $AB\iv$ admits a factorization of the above kind.  We remark
in this connection that
\bqn
S(A,B) &=& \Big(I+PM(AB\iv)Q\Big)\Big(T(AB\iv)+Q\Big) M(B),
\eqn
where $I+PM(AB\iv)Q$ and $M(B)$ are invertible operators.
Hence the problem of computing the dimension of the kernel and cokernel
of a singular integral operator can be reduced to a factorization problem
with the determination of partial indices.

Fredholm criteria related to a factorization problem and formulas for
the dimension of the kernel and cokernel similar to above have so
far not been known for singular integral operators with flip,
\begin{equation}\label{f1.SIOF}
PM(A)+PJM(B)+QJM(C)+QM(D),\qquad
A,B,C,D\in L^\iy(\T)\NN,
\end{equation}
not even in the case where the generating functions are smooth.
Also for Toeplitz + Hankel operators,
\begin{equation}
T(A)+H(B),\qquad A,B\in L^\iy(\T)\NN,
\end{equation}
such results have not yet been obtained.
Here 
\bqn
H(B) &=& PM(B)JP
\eqn
stands for the Hankel operator acting on $(H^2)^N$ with the generating
function $B\in L^\iy(\T)\NN$.
Only in a recent paper of E.~L.~Basor and the author \cite{BE99}
it has been observed that the invertibility of special class of
Toeplitz + Hankel operators might be related to a factorization problem.

The Fredholm theory of Toeplitz + Hankel operators with piecewise continuous
functions can be found in \cite{P} (see also \cite[Sect.~4.95--4.102]{BS89}).
Several aspects of the Fredholm theory of singular integral operators with flip
(also in a different settings) can be found in the monograph \cite{KL}.

We remark that there exists a ``classical'' trick, which allows
to reduce singular integral operators with flip to singular integral
operators without flip (and thus to a factorization problem).
This trick will be sketched below. Unfortunately, this trick
leads only to sufficient conditions and gives in general only
estimates on the dimensions of the kernel and cokernel.

The purpose of this paper is to consider general singular integral
operators with flip and Toeplitz + Hankel operators with sufficiently
smooth (e.g., H\"older continuous) matrix valued generating functions.
In the case where these operators are Fredholm we will establish
formulas for the dimension of the kernel and cokernel.  Note that (in
the case of continuous generating functions) Fredholm criteria are
easy to obtain.  These formulas will rely on a factorization problem,
which is slightly different from the classical Wiener--Hopf factorization.
Instead of
the partial indices $\ka_1,\dot,\ka_N\in\Z$, a collection of pairs
$(\ro_1,\ka_1),\dots(\ro_N,\ka_N)\in\{-1,1\}\times\Z$ appears, which
contain the relevant information about the dimension of the kernel and
cokernel and allow us to give an answer to the invertibility problem.

The general case, i.e., Fredholm criteria in terms of a factorization
problem for singular integral operators with flip and
Toeplitz + Hankel operators with generating functions in
$L^\iy(\T)\NN$, will be deferred to a future paper.

Let us state some basic relations between the operators introduced above.
Obviously, $P^2=P$, $Q^2=Q$ and $P+Q=I$ by definition. Moreover,
\begin{equation}
J^2=I,\qquad
JPJ=Q\quad\mbox{ and }\quad
JM(A)J=M(\wt{A}),
\end{equation}
where $\wt{A}$ stands for the function defined by
\bqn
\wt{A}(t) &=& A(1/t),\qquad t\in\T.
\eqn
For functions $A,B\in L^\iy(\T)\NN$, the following relation for multiplication
operators holds:
\bqn
M(AB) &=& M(A)M(B).
\eqn
From this and the above relations, one can deduce 
well known identities for Toeplitz and Hankel operators:
\bqn\label{f.Tab}
T(AB) &=& T(A)T(B)+H(A)H(\wt{B}),\\
\label{f.Hab}
H(AB) &=& T(A)H(B)+H(A)T(\wt{B}).
\eqn

Now let us explain how the above mentioned ``classical'' trick works
in regard to singular integral operators with flip. It works, of course,
also for Toeplitz + Hankel operators. First consider the identity
\bqn\label{f1.XY}
\frac{1}{2}
\left(\ba{cc} I & I \\ J & -J\ea\right)
\left(\ba{cc}X+YJ & 0 \\ 0 & X-YJ\ea\right)
\left(\ba{cc} I & J \\ I &-J\ea\right) &=&
\left(\ba{cc} X & Y \\ JYJ & JXJ\ea\right),
\eqn
where $X$ and $Y$ are arbitrary operators acting on $(L^2)^N$.
Note that the block operators on the left and the right of the
left hand side of the equation are the inverses of each other.
Given $a,b,c,d\in L^\iy(\T)\NN$, write
\bqn
A &=& \left(\ba{cc} a & b \\ c & d \ea\right)\in L^\iy(\T)\NNN
\eqn
and introduce two singular integral operators with flip:
\bqn\label{f1.Phi+}
\Phi(A) &=&
PM(a)+PM(b)J+QM(\ti{d})+QM(\ti{c})J,\\
\label{f1.Phi-}
\Phi'(A) &=&
PM(a)-PM(b)J+QM(\ti{d})-QM(\ti{c})J.
\eqn
Notice the slight change in notation in comparison with (\ref{f1.SIOF}).
With $X=PM(a)+QM(\ti{d})$ and $Y=PM(b)+QM(\ti{c})$ we can employ
(\ref{f1.XY}), and it follows that problem of Fredholmness,
invertibility and dimension of the kernel and cokernel are the same for
the operators
\bqn
\left(\ba{cc}\Phi_+(A) & 0 \\ 0 & \Phi_-(A)\ea\right)
&\mbox{ and }&
\left(\ba{cc} PM(a)+QM(\ti{d}) & PM(b)+QM(\ti{c}) \\
QM(\ti{b})+PM(c) & QM(\ti{a})+PM(d) \ea\right).
\eqn
However, this last operator can be rewritten as
\bqn\label{f1.SIO}
P\left(\ba{cc} M(a) & M(b) \\ M(c) & M(d) \ea\right)+
Q\left(\ba{cc} M(\ti{d}) & M(\ti{c}) \\
M(\ti{b}) & M(\ti{a}) \ea\right)
&=&
PM(A)+QM(W\wt{A}W)
\eqn
with a constant $2N\times 2N$ matrix
\bqn
W&=&\left(\ba{cc} 0 & I \\ I & 0 \ea\right).
\eqn
The operator (\ref{f1.SIO}) is a usual singular integral operator
with generating functions of twice the original matrix size.
By what has been said above about singular integral operators,
one is led to the factorization of matrix function $AW\wt{A}\iv W$
in the form (\ref{f1.fact}).

The disadvantage of this trick is that one cannot study
$\Phi(A)$ alone, but one is compelled to take also the ``conjugate''
operator $\Phi'(A)$ into account. In the worst case it can happen that
$\Phi(A)$ is a Fredholm operator whereas $\Phi'(A)$ is not, in which case
one obtains no information at all about $\Phi(A)$.

\section{First results about Toeplitz + Hankel operators}
\label{s2}

In this section we first establish the basic properties of general
Toeplitz + Hankel operators $T(A)+H(B)$ with $A,B\in L^\iy(\T)\NN$.
Then we introduce two special classes of such Toeplitz + Hankel operators
and consider their basic properties, too.  The further study of these
particular as well as of the general Toeplitz + Hankel operators will
be continued in later sections.

The following necessary condition for the Fredholmness of general Toeplitz +
Hankel operators is certainly well known.  For completeness sake, we present it
with a proof.
\begin{proposition}\label{p2.2}
Let $A,B\in L^\iy(\T)\NN$, and assume that $T(A)+H(B)$ is Fredholm.
Then $A\in G(L^\iy(\T)\NN).$
\end{proposition}
\begin{proof}
If $T(A)+H(B)$ is Fredholm, then there exist $\delta>0$ and a finite
rank projection $K$ on the kernel of $T(A)+H(B)$ such that 
\bqn
\|T(A)f+H(B)f\|_{(H^2)^N}+\|Kf\|_{(H^2)^N} &\geq& \delta\|f\|_{(H^2)^N}\nn
\eqn
for all $f\in(H^2)^N$. Replacing $f$ by $Pf$ and applying the estimate
$\|Pf\|\ge\|f\|-\|Qf\|$, it follows that
\bqn
\|T(A)f+H(B)f\|_{(L^2)^N}+\|KPf\|_{(L^2)^N}+\delta
\|Qf\|_{(L^2)^N} &\geq& \delta\|f\|_{(L^2)^N}.\nn
\eqn
for all $f\in(L^2)^N$. Introducing the isometries $U_{\pm n}:(L^2)^N\to
(L^2)^N,f(t)\mapsto t^{\pm n}f(t)$ and replacing $f$ by $U_nf$, we obtain
$$
\|U_{-n}T(A)U_nf+U_{-n}H(B)U_nf\|_{(L^2)^N}+\|KPU_nf\|_{(L^2)^N}
+\delta
\|U_{-n}QU_nf\|_{(L^2)^N} \geq \delta\|f\|_{(L^2)^N}.\nn
$$
Because $U_{\pm n}$ commute with multiplication operators and
$U_nJ=JU_{-n}$, we can write
$$
U_{-n}T(A)U_{n}=U_{-n}PU_nM(A)U_{-n}PU_n,\qquad
U_{-n}H(B)U_{n}=U_{-n}PU_{-n}M(B)JU_{-n}PU_n.
$$
Observe that $U_{-n}PU_n\to I$ and $U_{-n}PU_{-n}\to 0$ strongly as $n\to\iy$.
Hence it follows that $U_{-n}T(A)U_{n}\to M(A)$ and $U_{-n}H(B)U_{n}\to 0$
strongly as $n\to\iy$. Now we take the limit in the above norm
estimate. Since $U_n\to0$ weakly and $K$ is compact, we have $KPU_n \to 0$
strongly. Moreover, $U_{-n}QU_n\to 0$ strongly. We obtain that
\bqn
\|M(A)f\|_{(L^2)^N} &\ge& \delta\|f\|_{(L^2)^N}\nn
\eqn
for all $f\in(L^2)^N$.
From this it immediately follows that $A\in G(L^\iy(\T)\NN)$.
\end{proof}

For continuous matrix valued functions $A$ and $B$, the just stated
necessary Fredholm condition is also sufficient.
Recall in this connection that the winding number of a complex valued
nonvanishing continuous functions $a$ defined on the unit circle
is given by
\bqn
\wind a &=& \Big[\frac{1}{2\pi}\arg a(e^{i\theta})\Big]_{\theta=0}^{2\pi},
\eqn
where the argument $\arg a(e^{i\theta})$ is chosen continuously
on $[0,2\pi]$. Again, the following result is well known, and we present the
proof only for completeness sake.

\begin{proposition}\label{p2.4}
Let $A,B\in C(\T)\NN$. Then $T(A)+H(B)$ is Fredholm if and only if
$A\in G(C(\T)\NN)$. Moreover, if this is true, then
$\ind(T(A)+H(B))=-\wind\det A$.
\end{proposition}
\begin{proof}
It suffices to remark that the Hankel operator with a continuous generating
function is compact. Hence, by making use of (\ref{f.Tab}), it is easy to see
that a Fredholm regularizer for $T(A)+H(B)$ is given by $T(A\iv)$.
As to the index formula, we remark that for $B\in C(\T)\NN$, $A\in G(C(\T)\NN)$,
$$
\ind(T(A)+H(B)) \;=\; \ind T(A) \;=\; \ind T(\det A) \;=\;
-\wind \det A.
$$
The last equality is the well known formula for the Fredholm index of
a scalar Toeplitz operator with continuous symbol.  For the precise
justification of the second last equality see, e.g.,
\cite[Thm.~2.94]{BS89} 
\end{proof}

After these results for general Toeplitz + Hankel operators we are
going to consider two special classes of Toeplitz + Hankel operators.
These operators possess a number of unexpected properties.

In what follows, let $W\in\C\NN$ be any matrix such that $W^2=I$.
For $A\in L^\iy(\T)\NN$, we introduce the operators
\bqn
\cM_W(A) &=& T(A)+H(AW),\label{f2.MW}\\
\cN_W(A) &=& T(A)+H(W\wt{A}).\label{f2.NW}
\eqn
What makes these classes of operators so interesting for us is the fact
that an analogue of formula (\ref{f.Tab}) holds. Indeed,
\bqn
\label{f.Mab}
\cM_W(AB) &=& \cM_W(A)\cM_W(B)+H(AW)\cM_W(W\wt{B}W-B),\\
\label{f.Nab}
\cN_W(AB) &=& \cN_W(A)\cN_W(B)+\cN_W(W\wt{A}W-A)H(W\wt{B}).
\eqn
These formulas can be verified straightforwardly by using
(\ref{f.Tab}), (\ref{f.Hab}), and the assumption that $W$ is a
constant matrix with $W^2=I$:
\bqn
\cM_W(AB) &=& T(AB)+H(ABW)\nn\\
&=& T(A)T(B)+H(A)H(\wt{B})+ T(A)H(BW)+ H(A)T(\wt{B}W)\nn\\
&=& T(A)\cM_W(B)+ H(A)\cM_W(\wt{B}W)\nn\\
&=& T(A)\cM_W(B)+ H(AW)\cM_W(W\wt{B}W)\nn\\
&=& \cM_W(A)\cM_W(B)+ H(AW)\cM_W(W\wt{B}W-B).\nn
\eqn
Similarly,
\bqn
\cN_W(AB) &=& T(AB)+H(W\wt{A}\wt{B})\nn\\
&=& T(A)T(B)+H(A)H(\wt{B})+ T(W\wt{A})H(\wt{B})+ H(W\wt{A})T(B)\nn\\
&=& \cN_W(A)T(B)+ \cN_W(W\wt{A})H(\wt{B})\nn\\
&=& \cN_W(A)T(B)+ \cN_W(W\wt{A}W)H(W\wt{B})\nn\\
&=& \cN_W(A)\cN_W(B)+\cN_W(W\wt{A}W-A)H(W\wt{B}).\nn
\eqn

Next we introduce the set
\bqn
(L^\iy)\NN_W &=& \Big\{\;A\in L^\iy(\T)\NN\;:\;W\wt{A}W=A\;\Big\}.
\eqn
We remark that $(L^\iy)\NN_W$ is an inverse closed Banach subalgebra
of $L^\iy(\T)\NN$.

Under additional assumptions on the functions $A$ or $B$, formulas
(\ref{f.Mab}) and (\ref{f.Nab}) can be simplied. Indeed,
\bqn
\cM_W(AB) &=&\lefteqn{\cM_W(A)\cM_W(B)}\hspace{23ex}
\mbox{if $A\in (\ovl{H^\iy})\NN$ or $B\in (L^\iy)\NN_W$};
\label{f2.Mmult}\\
\cN_W(AB) &=&\lefteqn{\cN_W(A)\cN_W(B)}\hspace{23ex}
\mbox{if $A\in (L^\iy)\NN_W$ or $B\in (H^\iy)\NN$}.
\label{f2.Nmult}
\eqn
Consequently, in some cases the mappings $A\mapsto\cM_W(A)$ and
$A\mapsto\cN_W(A)$ are multiplicative. This is halfway not surprising.
In fact,
\bqn
\cM_W(A)&=&T(A)\qquad
\mbox{if $A\in (\ovl{H^\iy})\NN$};
\label{f2.M=T}
\\
\cN_W(A)&=&T(A)\qquad
\mbox{if $A\in (H^\iy)\NN$}.
\label{f2.N=T}
\eqn
Hence in these cases we are dealing just with usual Toeplitz
operators which have symbols in $(\ovl{H^\iy})\NN$ and $(H^\iy)\NN$,
respectively.

More interesting is the case where $A\in(L^\iy)_W\NN$. It turns out that
then both of the above types of operators coincide:
\bqn\label{f2.M=N}
\cM_W(A) &=& \cN_W(A)\qquad\mbox{ if $A\in (L^\iy)_W\NN$}.
\eqn
Moreover, the following result shows that both the invertibility
and the Fredholm problem can be solved completely in a very simple way.

\begin{corollary}\label{c2.3}
Let $A\in (L^\iy)_W\NN$. Then the following is equivalent:
\begin{itemize}
\item[(i)]
$A\in G(L^\iy)_W\NN$;
\item[(ii)]
$\cM_W(A)=\cN_W(A)$ is invertible;
\item[(iii)]
$\cM_W(A)=\cN_W(A)$ is Fredholm.
\end{itemize}
If this is fulfilled, then the inverse of $\cM_W(A)=\cN_W(A)$
is given by $\cM_W(A\iv)=\cN_W(A\iv)$.
\end{corollary}
\begin{proof}
Because of the multiplicative relations (\ref{f2.Mmult}) or
(\ref{f2.Nmult}), it follows that (i) implies (ii), where the inverse
of $\cM_W\iv(A)=\cN_W\iv(A)$ is given by $\cM_W(A\iv)=\cN_W(A\iv)$.
The implication (ii)$\Rightarrow$(iii) is obvious.  The fact that
(iii) implies (i) follows from Proposition \ref{p2.2} in connection
with the inverse closedness of $(L^\iy)_W\NN$ in $L^\iy(\T)\NN$.
\end{proof}

The operators $\cM_W(A)$ and $\cN_W(A)$ are not completely unrelated
with each other. First of all, there is a connection by means of the
adjoints,
\begin{equation}
(\cM_W(A))^* \;=\; \cN_{W^*}(A^*)\quad\mbox{ or }\quad
(\cN_W(A))^* \;=\; \cM_{W^*}(A^*),
\end{equation}
where $A^*(t):=(A(t))^*$. Here we need only remark that
$P^*=P$, $J^*=J$ and $M(A)^*=M(A^*)$, from which
$T(A)^*=T(A^*)$ and $H(A)^*=H(\wt{A}^*)$ follows.

Another relation is established by the identity
\bqn\label{f2.MWNW}
\cM_W(A)\cN_W(B) &=& T(AB)+H(AW\wt{B})
\eqn
Indeed,
\bqn
\cM_W(A)\cN_W(B) &=&
\Big(T(A)+H(AW)\Big)\Big(T(B)+H(W\wt{B})\Big)\nn\\
&=& T(A)T(B)+H(A)H(\wt{B})+H(AW)T(B)+T(AW)H(\wt{B})\nn\\
&=& T(AB)+H(AW\wt{B}).\nn
\eqn
Here we have only used the assumption that $W$ is a constant matrix
with $W^2=I$ and formulas (\ref{f.Tab}) and (\ref{f.Hab}).

Finally, we illustrate some further interesting consequences of the relations
(\ref{f2.Mmult}) and (\ref{f2.Nmult}).

\begin{corollary}
Let $A\in L^\iy(\T)\NN$.
\begin{itemize}
\item[(i)]
If $A$ admits a factorization $A(t)=A_-(t)A_0(t)$ with
$A_-\in G(\ovl{H^\iy})\NN$ and $A_0\in G(L^\iy)_W\NN$, then $\cM_W(A)$
is invertible and the inverse equals $\cM_W(A_0\iv)T(A_-\iv)$.
\item[(ii)]
If $A$ admits a factorization $A(t)=A_0(t)A_+(t)$ with
$A_0\in G(L^\iy)_W\NN$ and $A_+\in G(H^\iy)\NN$, then $\cN_W(A)$
is invertible and the inverse equals $T(A_+\iv)\cN_W(A_0\iv)$.
\end{itemize}
\end{corollary}
\begin{proof}
As to assertion (i), it follows from (\ref{f2.Mmult}) and (\ref{f2.M=T}) that
$\cM_W(A)=T(A_-)\cM_W(A_0)$. The inverse of $T(A_-)$ equals $T(A_-\iv)$ and
the inverse of $\cM_W(A_0)$ equals $\cM_W(A_0\iv)$. In regard to assertion (ii),
we use (\ref{f2.Nmult}) and (\ref{f2.N=T}) and  obtain
$\cN_W(A)=T(A_+)\cN_W(A_0)$. The inverse of $T(A_-)$ equals $T(A_-\iv)$
and the inverse of $\cN_W(A)$ equals $\cN_W(A_0\iv)$.
\end{proof}

We conclude this section by making some more or less heuristic
remarks, which should serve as a motivation for the kinds of
factorizations that we are going to consider in the following section.

The above necessary condition for the invertibility of $\cM_W(A)$ (and likewise
for $\cN_W(A)$) is certainly for several reason far away from being sufficient.
In analogy to the usual theory of the Wiener--Hopf factorization one may
guess that under certain conditions there exist a factorization of the form
\bqn\label{f2.2facts}
A(t) &=& A_-(t)R(t)A_0(t)
\eqn
with appropriate conditions on the factors $A_-$, $A_0$,
and where the middle factor $R$ is of a particularly simple form.
Indeed, if $A_-\in G(\ovl{H^\iy})\NN$ and $A_0\in G(L^\iy)_W\NN$,
then the invertibility of $\cM_W(A)$ is equivalent to the
invertibility of $\cM_W(R)$. More general, the dimensions of the
kernel and cokernel of $\cM_W(A)$ coincide with those for
$\cM_W(R)$. If $R$ is of a particularly simple form, then one can hope
that these dimensions can be calculated.

It turns out that the factorization of the form (\ref{f2.2facts}),
which may deserve the name ``asymmetric'', can be related
to some kind of Wiener--Hopf factorization, which looks kind of
``antisymmetric''. Indeed, if we are given the factorization
(\ref{f2.2facts}), then
\begin{equation}
AW\wt{A}\iv \;=\;
A_-R(t)A_0W\wt{A}_0\iv \wt{R}\iv\wt{A}_-\iv
\;=\; A_-RW\wt{R}\iv\wt{A}_-\iv,
\end{equation}
where the last equality follows from the presumed property $A_0=W\wt{A}_0W$
of the factor $A_0$. Replacing the product $R(t)W\wt{R}\iv(t)$ by the notation
$D(t)$, we arrive at a factorization 
\bqn\label{f2.anti-fac}
F(t) &=& A_-(t)D(t)\wt{A}_-\iv(t)
\eqn
of the matrix function
\bqn
F(t) &=& A(t)W\wt{A}\iv(t).
\eqn
Assuming for a moment that $D(t)$ is of an appropriate form,
it follows that (\ref{f2.anti-fac}) is some kind of Wiener--Hopf factorization,
where the right and the left factors are related with each other in an 
``antisymmetric'' way.

Similarly, the analysis of the operator $\cN_W(A)$ may lead to
a factorization of the form
\bqn\label{f2.2facts2}
A(t) &=& A_0(t)R(t)A_+(t),
\eqn
again with suitable conditions on the factors. Elaborating on this
``asymmetric'' factorization, we arrive at the following
``antisymmetric'' factorization,
\bqn\label{f2.anti-fac2}
G(t) &=& \wt{A}_+\iv(t)D(t)A_+(t)
\eqn
of the matrix function
\bqn
G(t) &=& \wt{A}\iv(t)WA(t).
\eqn
Here $D(t)$ stands for $\wt{R}\iv(R)WR(t)$, which slightly differs from
the previous situation.

The reader should observe that whereas the ``asymmetric'' factorizations
(\ref{f2.2facts}) and (\ref{f2.2facts2}) are of different types,
the ``antisymmetric'' factorizations (\ref{f2.anti-fac}) and
(\ref{f2.anti-fac2}) is essentially of the same form
only that different notation has been used.

\section{Some results about factorizations}\label{s3}

\subsection{The usual factorization within a Banach algebra}

Throughout the rest of this paper let $\cB$ stand for a Banach algebra of
functions defined on the unit circle such that the following properties
are fulfilled:
\begin{itemize}
\item[(a)]
$\cB$ is an inverse closed Banach subalgebra of $C(\T)$;
\item[(b)]
$\cB$ contains all trigonometric polynomials;
\item[(c)]
If $a\in\cB$, then $\wt{a}\in\cB$;
\item[(d)]
For each $N$, each matrix function $A\in G\cB\NN$ admits
a factorization of the form 
\bqn\label{f3.Afact}
A(t) &=& A_-(t)\Lambda(t)A_+(t)
\eqn
where $\Lambda(t)=\diag(t^{\ka_1},t^{\ka_2},\dots,t^{\ka_N})$ with
$\ka_1,\dots,\ka_N\in\Z$,
\begin{equation}\label{f3.Apm}
A_+\in G\cB_+\NN\quad\mbox{ and }\quad A_-\in G\cB_-\NN,
\end{equation}
where
\begin{equation}\label{f3.Bpm}
\cB_+\;:=\;\cB\cap H^\iy\quad\mbox{ and }\quad
\cB_-\;:=\;\cB\cap \ovl{H^\iy}.
\end{equation}
\end{itemize}

Examples of Banach algebras $\cB$ having the properties (a)--(d) are
the Wiener algebra $W$ or the Banach algebras $C^\alpha$ of all H\"older
continuous functions defined on the unit circle with exponents $0<\alpha<1$.

From the factorization point of view, only the assumptions (b) and (d)
and the condition that $\cB$ is a Banach subalgebra of $L^\iy(\T)$ are
important.  More specifically, one refers to the factorization
(\ref{f3.Afact}) with the properties (\ref{f3.Apm}) as a factorization
within the Banach algebra $\cB$.  Related to this concept are such
notions as that of decomposing Banach algebras and Banach algebras
with factorization property.  We will go into these details, but
simply refer the reader to \cite[Sect.~10.14--10.23]{BS89}.  We also
note that $\cB_+$ and $\cB_-$ defined in (\ref{f3.Bpm}) are Banach
subalgebras of $\cB$ containing the unit element.

It is obvious that a factorization in such a Banach algebra is
automatically a generalized factorization (or, $\Phi$-factorization)
in the space $L^2$.  In particular, $A_+\in G(H^\iy)\NN$ and $A_-\in
G(\ovl{H^\iy})\NN$, and thus the factors $A_+$ and $A_-$ satisfy the
conditions (i)--(iii) stated in the introduction.

Our assumption (a) is motivated by the circumstance that we will confine
ourselves to continuous matrix valued functions because in this case
Fredholm criteria for Toeplitz + Hankel operators and singular integral
operators with flip are easy to obtain. The inverse closedness is
needed for the conclusion that each function $A\in\cB\NN$ which is
invertible (in $C(\T)\NN$) admits a factorization of the above kind.

The assumption (c) will be important for the definition of another
type of factorization that we will introduced later on.  We remark
in this connection the obvious fact that $A\in\cB_+\NN$ if and only if
$\wt{A}\in\cB_-\NN$. Consequently, $A\in G\cB_+\NN$ if and only if
$\wt{A}\in G\cB_-\NN$.

As has already been noted in the introduction, the partial indices of
such factorizations are uniquely determined up to change of order.  In
fact, the order of the partial indices can be changed in any desired
way.  Namely, one can replace $F_-(t)$ with $F_-(t)\Pi\iv$,
$\Lambda(t)$ with $\Pi\Lambda(t)\Pi\iv$ and $F_+(t)$ with $\Pi
F_+(t)$, where $\Pi$ is a suitable permutation matrix.

The following result is well known \cite{CG,LS} and answers the
question about the uniqueness of the factors $A_+$ and $A_-$ in a
factorization.  In order to simplify the statement we will assumes
without loss of generality that the partial indices are ordered
increasingly.  Then the factors corresponding to different
factorizations are related with each other by certain rational block
triangular matrix functions whose structure is determined by the
multiple occurrence of same values for the partial indices.  In this
regard, we introduce the notation $I_l$ for the identity matrix of size
$l\times l$.

\begin{proposition}\label{p3.1}
Assume that we are given two factorizations of a function $F\in G\cB\NN$,
\begin{equation}
F(t) \;=\; F_-^{(1)}(t)\Lambda(t)F_+^{(1)}(t)
\;=\; F_-^{(2)}(t)\Lambda(t)F_+^{(2)}(t)
\end{equation}
with $F_-^{(j)}\in G\cB_-\NN$, $F_+^{(j)}\in G\cB_+\NN$, and
\begin{equation}\label{f3.La}
\Lambda(t)\;=\;
\diag(t^{\bka_1}I_{l_1},t^{\bka_2}I_{l_2},\dots,t^{\bka_R}I_{l_R}),
\end{equation}
where $R\in\{1,2,\dots\}$, $l_1,\dots,l_R\in\{1,2,\dots\}$,
$l_1+\dots+l_R=N$, $\bka_1,\dots,\bka_R\in\Z$ and
\begin{equation}
\bka_1\;<\;\bka_2\;<\dots <\;
\bka_{N-1}\;<\;\bka_R.
\end{equation}
Then there exist matrix functions $U$ and $V$ which are of the form
\begin{equation}\label{f3.UV}
U(t)=\left(
\ba{cccc} A_{11} & U_{12}(t) & \cdots & U_{1R}(t) \\
0 & A_{22} &\ddots& \vdots \\ \vdots &\ddots& \ddots & U_{R-1,R}(t)\\
0&\cdots&0& A_{RR}\ea\right),\;\;
V(t)=\left(
\ba{cccc} A_{11} & V_{12}(t) & \cdots & V_{1R}(t) \\
0 & A_{22} &\ddots& \vdots \\ \vdots &\ddots& \ddots &V_{R-1,R}(t) \\
0&\cdots&0& A_{RR}\ea\right)
\end{equation}
with $A_{jj}\in G\C^{l_j\times l_j}$ and
\begin{equation}
U_{jk}(t) \;=\; \sum_{m=0}^{\bka_k-\bka_j}A_{jk}^{(m)}t^m,\quad
V_{jk}(t) \;=\; t^{\bka_j-\bka_k}U_{jk}(t),\quad
A_{jk}^{(m)}\in\C^{l_j\times l_k}
\end{equation}
for $1\le j<k\le R$ such that
\begin{equation}\label{f3.Fpm12}
F_-^{(2)}(t)\;=\;F_-^{(1)}(t)V(t),\qquad
F_+^{(1)}(t)\;=\;U(t)F_+^{(2)}(t).
\end{equation}
\end{proposition}

Due to the assumption (b) on $\cB$, it is not hard too see that $U\in
G\cB_+\NN$ and $V\in G\cB_-\NN$.  The previous proposition holds, by
the way, not only for factorizations within the Banach algebra $\cB$,
but also for generalized factorizations (see (\ref{f1.fact})
and (i)--(iii)).  However, we will not make
use of this fact.

Actually, the statement of this proposition can be reversed.
If we are given a factorization $F(t)=F_-^{(1)}(t)\Lambda(t)F_+^{(1)}(t)$,
introduce functions $U$ and $V$ with the above properties and define
$F_-^{(2)}$ and $F_+^{(2)}$ by (\ref{f3.Fpm12}), then
also $F(t)=F_-^{(2)}(t)\Lambda(t)F_+^{(2)}(t)$ is such a factorization.

\subsection{Antisymmetric factorization within a Banach algebra}

In what follows we are going to introduce and study a slightly different type
of factorization. It is essentially also a factorization of the form
(\ref{f3.Afact}), but we require in addition that the factors $F_+$ and $F_-$
are related with each other by $F_+(t)=\wt{F}_-\iv(t)$.
Moreover, the middle factor is allowed to be of a more general form. Namely,
\bqn\label{f3.D(t)}
D(t) &=& \diag(\ro_1t^{\ka_1},\ro_2t^{\ka_2},\dots,\ro_Nt^{\ka_N})
\eqn
with $\ka_1,\dots,\ka_N\in\Z$ and $\ro_1,\dots,\ro_N\in\{-1,1\}$.

More specifically, we are going to consider a factorization of a function
$F\in G\cB\NN$ in the form
\begin{equation}\label{f3.fac2}
F(t)\;=\;F_-(t)D(t)\wt{F}_-\iv(t)\quad\mbox{ with }\quad
F_-\in G\cB_-\NN,
\end{equation}
where $D(t)$ is given by (\ref{f3.D(t)}). Such a factorization will be
called an {\em antisymmetric factorization of $F$ within the Banach
algebra $\cB$}.

It will turn out that the collection of pairs
\begin{equation}\label{f3.pairs1}
(\ro_1,\ka_1),\;(\ro_2,\ka_2),\;\dots\;
(\ro_N,\ka_N)\;\in\;\{-1,1\}\times\Z
\end{equation}
plays the same important role as the collection of the partial indices
$\ka_1,\dots,\ka_N\in\Z$ in the classical situation. Therefore, we will call
this collection the {\em characteristic pairs} of the antisymmetric
factorization of $F$.

We first study the existence of an antisymmetric factorization.
Because $\wt{D}\iv(t)=D(t)$ for each middle factors of the above kind,
it is easy to see that the condition
\bqn
F(t) &=& \wt{F}\iv(t)
\eqn
is necessary for the existence of an antisymmetric factorization
of a function $F$. The following theorem shows that, essentially, this
condition is also sufficient.

\begin{theorem}\label{t3.2}
Assume that $F\in G\cB\NN$ satisfies the condition $F(t)=\wt{F}\iv(t)$.
Then there exists a function $F_-\in G\cB_-\NN$ such that $F$ can be
factored in the form
\bqn\label{f2.fac}
F(t) &=& F_-(t)D(t)\wt{F}_-\iv(t),
\eqn
where $D(t)$ is given by (\ref{f3.D(t)}) with certain characteristic pairs
(\ref{f3.pairs1}).
\end{theorem}
\begin{proof}
Because of the assumptions on the Banach algebra $\cB$ there exists
a factorization
\bqn\label{f2.facF1}
F(t) &=& F_-(t)\Lambda(t)F_+(t)
\eqn
with $F_\pm\in G\cB\NN_\pm$, $\Lambda(t)=\diag(t^{\ka_1},
\dots,t^{\ka_N})$, $\ka_1,\dots,\ka_N\in\Z$ and
$\ka_1\le\dots\le\ka_N$ without loss of generality. Taking the inverse and
replacing $t$ by $1/t$, it follows that
\bqn\label{f2.facF2}
\wt{F}\iv(t) &=& \wt{F}_+\iv(t)\Lambda(t)\wt{F}_-\iv(t).
\eqn
Because $\wt{F}\iv=F$, the expressions (\ref{f2.facF1}) and
(\ref{f2.facF2}) are equal and represent factorizations of the form
(\ref{f3.Afact}). We can apply Proposition \ref{p3.1} and write
$\Lambda(t)$ in the form (\ref{f3.La}) with the conditions on parameters
$R$, $l_1,\dots,l_R$ and $\bka_1,\dots,\bka_N$ stated there.
We conclude that there exists a matrix function $U(t)$ of
the form (\ref{f3.UV}) such that
\bqn\label{f2.FS}
F_+(t) &=& U(t)\wt{F}_-\iv(t),
\eqn
Combining (\ref{f2.FS}) with (\ref{f2.facF1}) and introducing
$X(t)=\Lambda(t)U(t)$, it follows that
\bqn\label{f2.facX}
F(t) &=& F_-(t)X(t)\wt{F}_-\iv(t),
\eqn
where $X(t)$ is of the form
\bqn\label{f2.Xmat}
X(t) &=& \left(
\ba{cccc} X_{11}(t) & X_{12}(t) & \cdots & X_{1R}(t) \\
0 & X_{22}(t) &\ddots&\vdots \\ \vdots &\ddots& \ddots & X_{R-1,R}(t)\\
0&\cdots&0& X_{RR}(t)\ea\right)
\eqn
with
$$
\ba{rclclcl}
X_{jk}(t) &=& \displaystyle\sum_{m=\bka_j}^{\bka_k}X_{jk}^{(m)}t^m, &&
X_{jk}^{(m)}\in\C^{l_j\times l_k},&&
1\le j<k\le R,\\[4ex]
X_{jj}(t) &=& X_jt^{\bka_j}, &&
X_{j}\in G\C^{l_j\times l_j},&&
1\le j\le R.\ea
$$
Introduce
\bqn\label{f2.X0}
X_0(t) &=& \diag(X_{11}(t),X_{22}(t),\dots,X_{RR}(t)),
\eqn
write
\begin{equation}\label{f2.XN}
X(t)\;=\; (I+N_1(t))X_0(t) \;=\;X_0(t)(I+N_2(t)),
\end{equation}
and observe that $N_1(t)$ and $N_2(t)$ are nilpotent matrix functions.
Note also that $N_1\in\cB_-\NN$ and $N_2\in\cB_+\NN$.  From formula
(\ref{f2.XN}) we obtain that $N_1(t)X_0(t)=X_0(t)N_2(t)$, and moreover
$N_1^m(t)X_0(t)=X_0(t)N_2^m(t)$ for each $m$ by induction.  The matrix
functions $(I+N_1(t))^{1/2}$ and $(I+N_2(t))^{1/2}$ are well defined
by a series expansion, which is finite due to the nilpotency.  Using this
series expansion, it follows that
\bqn
(I+N_1(t))^{1/2}X_0(t) &=& X_0(t)(I+N_2(t))^{1/2}.\nn
\eqn
This in connection with (\ref{f2.XN}) implies 
\bqn\label{f2.XNN}
X(t) &=& (I+N_1(t))^{1/2}X_0(t)(I+N_2(t))^{1/2}.
\eqn
From (\ref{f2.facX}) and the assumption $\wt{F}\iv(t)=F(t)$, it
follows that $\wt{X}\iv(t)=X(t)$. From the representation
(\ref{f2.Xmat}) and the definition (\ref{f2.X0}), we further obtain
$\wt{X}_0\iv(t)=X_0(t)$. On account of (\ref{f2.XN}), it now follows that
\bqn
I=\wt{X}(t)X(t) = (I+\wt{N}_1(t))\wt{X}_0(t)X_0(t)(I+N_2(t))
=(I+\wt{N}_1(t))(I+N_2(t)).\nn
\eqn
Hence $I+N_2(t)=(I+\wt{N}_1(t))\iv$, and, consequently,
$(I+N_2(t))^{1/2}=(I+\wt{N}_1(t))^{-1/2}$. This in connection with
(\ref{f2.XNN}) implies that
\bqn
\label{f2.XNN2}
X(t) &=& (I+N_1(t))^{1/2}X_0(t)(I+\wt{N}_1(t))^{-1/2}.
\eqn
From $\wt{X}_0\iv(t)=X_0(t)$ it follows (by putting $t=1$ and $t=-1$) that
$$
(X_0(1))^2\;=\;(X_0(-1))^2\;=\;I.
$$
Hence $X_j^2=I$ for each $1\le j\le R$. Thus we can write
$X_j=T_j\,\diag(I_{p_j},-I_{q_j})\,T_j\iv$ with certain
$T_j\in G\C^{l_j\times l_j}$ and $p_j,q_j\in\{0,1,\dots\}$ such that
$p_j+q_j=l_j$. It follows that
\bqn\label{f3.xxx}
X_0(t) &=& T\;\diag(\ro_1 t^{\ka_1},\ro_2 t^{\ka_2},\dots,\ro_N t^{\ka_N})\;T\iv
\eqn
with certain $\ro_1,\dots,\ro_N\in\{-1,1\}$ where
$T=\diag(T_1,T_2,\dots,T_R)\in G\C\NN$. Denoting by $D(t)$ the
diagonal matrix in (\ref{f3.xxx}), it follows in connection
with (\ref{f2.facX}) and (\ref{f2.XNN2}) that
\bqn
F(t) &=& F_-(t)(I+N_1(t))^{1/2}\,T\, D(t)\, T\iv\,
(I+\wt{N}_1(t))^{-1/2}\wt{F}_-\iv(t).\nn
\eqn
Because $(I+N_1(t))^{1/2}\in G\cB_-\NN$, we may replace the expression
$F_-(t)(I+N_1(t))^{1/2}\,T$ by the notation $F_-(t)$.
In this way, we arrive at the desired factorization (\ref{f2.fac}). 
\end{proof}

We remark that an antisymmetric factorization (\ref{f2.fac}) is obviously
an antisymmetric factorization of the form
\bqn
F(t) &=& \wt{F}_+\iv D(t) F_+(t)
\eqn
with $F_+\in G\cB_+\NN$ and the same middle factor $D(t)$. The only
difference is that of the notation of the factors.
Indeed, $F_+(t)=\wt{F}_-\iv(t)$ shows the relation.

The next theorem concerns the uniqueness of the characteristic
pairs of an antisymmetric factorization up to change of order.
Notice first that it is possible to rearrange the order
of these pairs in any desired way.
Indeed, one can replace $F_-(t)$ with $F_-(t)\Pi\iv$ and
$D(t)$ with $\Pi D(t)\Pi\iv$, where $\Pi$ is a suitable permutation matrix.
The important point is that the replacement of
$F_-(t)$ with $F_-(t)\Pi\iv$ implies the replacement of
$\wt{F}_-\iv(t)$ with $\Pi\wt{F}_-\iv(t)$, which fits with the
factorization formula (\ref{f3.fac2}).

\begin{theorem}\label{t3.3}
In an antisymmetric factorization of a function $F\in G\cB\NN$,
the characteristic pairs are uniquely determined up to change of order.
\end{theorem}
\begin{proof}
Because an antisymmetric factorization is automatically
also a usual factorization in the Banach algebra $\cB$ (except for the
slightly different middle factor, which is irrelevant at this place),
it follows that the numbers $\ka_1,\dots,\ka_N$ are uniquely determined up
to change of order. Because the order of the characteristic pairs
in an antisymmetric factorization can be rearranged in any desired way,
we can assume without loss of generality that $\ka_1\le\dots\le\ka_N$.

Now suppose that we are given two antisymmetric factorizations of $F$,
$$
F(t) \;=\; F_-^{(1)}(t) D^{(1)}(t) (\wt{F}_-^{(1)})\iv(t)
\;=\; F_-^{(2)}(t) D^{(2)}(t) (\wt{F}_-^{(2)})\iv(t),
$$
where $D^{(1)}(t)$ and $D^{(2)}(t)$ are both of the form (\ref{f3.D(t)})
but with the pairs
$$
(\ro_1^{(1)},\ka_1),\dots,(\ro_N^{(1)},\ka_N)\quad\mbox{ and }\quad
(\ro_1^{(2)},\ka_1),\dots,(\ro_N^{(2)},\ka_N),
$$
respectively. Introducing the parameters $l_1,\dots,l_R\in\{1,2,\dots\}$
and the integers $\bka_1,\dots,\bka_R$ as in Proposition \ref{p3.1},
we can write 
\bqn
D^{(1)}(t) &=& \diag(S_1^{(1)}t^{\bka_1},S_2^{(1)}t^{\bka_2},\dots,
S_R^{(1)}t^{\bka_R}),\nn\\
D^{(2)}(t) &=& \diag(S_1^{(2)}t^{\bka_1},S_2^{(2)}t^{\bka_2},\dots,
S_R^{(2)}t^{\bka_R}),\nn
\eqn
where $S_k^{(1)}$ and $S_k^{(2)}$ are diagonal matrices of size
$l_k\times l_k$ with entries $-1$ or $1$ on the diagonal. 
Moreover, we can write $D^{(1)}(t)=\Lambda(t)S^{(1)}$ and 
$D^{(2)}(t)=\Lambda(t)S^{(2)}$,
where $\Lambda(t)$ is of the form (\ref{f3.La}) and
$S^{(j)}=\diag(S^{(j)}_1,S^{(j)}_2,\dots,S^{(j)}_R)$.
It follows that
$$
F(t)\;=\; F_-^{(1)}(t) \Lambda(t)\Big(S^{(1)}(\wt{F}_-^{(1)})\iv(t)\Big)
\;=\; F_-^{(2)}(t) \Lambda(t)\Big(S^{(2)}(\wt{F}_-^{(2)})\iv(t)\Big),
$$
are two factorizations of the form (\ref{f3.Afact}). We apply Proposition
\ref{p3.1} and see that
$$
F_-^{(2)}(t)\;=\;F_-^{(1)}(t)V(t),\qquad
S^{(1)}(\wt{F}_-^{(1)})\iv(t)\;=\;U(t)S^{(2)}(\wt{F}_-^{(2)})\iv(t),
$$
where $U$ and $V$ are of the form (\ref{f3.UV}).
The last equation can be rewritten as $F_-^{(1)}(t)(S^{(1)})\iv
=F_-^{(2)}(t)(S^{(2)})\iv\wt{U}\iv(t)$. Combined with the first equation,
it follows that
$$
S^{(1)}\;=\;\wt{U}(t)S^{(2)}V\iv(t).
$$
Because of the block triangular structure of $U$ and $V$ with invertible
constant matrices $A_k$ on the block diagonal, we obtain that
$S_k^{(1)}=A_kS^{(2)}A_k\iv$ for each $k=1,\dots,R$. Hence $S_k^{(1)}\sim
S_k^{(2)}$, and, consequently, the numbers of $1$'s and $-1$'s, respectively,
on the diagonal of $S_k^{(1)}$ and $S_k^{(2)}$ is the same. 
From this it follows that the collection of the pairs $(\ro_k^{(1)},\ka_k)$
is the same as the collection of the pairs $(\ro_k^{(2)},\ka_k)$
up to change of order.
\end{proof}

It is possible (similar as has been done in Proposition \ref{p3.1}) to
state the relation between the factors $F_-$ of two different antisymmetric
factorizations of a given function. We will omit this result because it
is a little bit difficult to state and will not be needed for
our purposes.

For a given antisymmetric factorization of a function $F$ with
characteristic pairs (\ref{f3.pairs1}), we introduce the
following nonnegative integers:
\begin{itemize}
\item[] $\alpha=$
number of $k\in\{1,\dots,N\}$ for which $\ro_k=1$ and $\ka_k$ is even;
\item[] $\beta=$
number of $k\in\{1,\dots,N\}$ for which $\ro_k=1$ and $\ka_k$ is odd;
\item[] $\gamma=$
number of $k\in\{1,\dots,N\}$ for which $\ro_k=-1$ and $\ka_k$ is odd;
\item[] $\delta=$
number of $k\in\{1,\dots,N\}$ for which $\ro_k=-1$ and $\ka_k$ is even.
\end{itemize}
Besides the obvious fact that $\alpha+\beta+\gamma+\delta=N$,
the following ``a priori'' characterization of these numbers can be
obtained.

\begin{proposition}\label{p3.4}
Assume that $F\in G\cB\NN$ admits an antisymmetric factorization with
the numbers $\alpha,\beta,\gamma,\delta$ be defined as above. Then 
\begin{equation}
F(1)\;\sim\;\diag(I_{\alpha+\beta},-I_{\gamma+\delta})\quad\mbox{ and }
\quad F(-1)\;\sim\;\diag(I_{\alpha+\gamma},-I_{\beta+\delta}).
\end{equation}
\end{proposition}
\begin{proof}
Putting $t=1$ or $t=-1$ in the factorization
$F(t)=F_-(t)D(t)\wt{F}_-\iv(t)$ it follows that
$F(1)\sim D(1)$ and $F(-1)\sim D(-1)$. Now the assertion follows from
the facts that $D(1)\sim\diag(I_{\alpha+\beta},-I_{\gamma+\delta})$ and
$D(-1)\sim\diag(I_{\alpha+\gamma},-I_{\beta+\delta})$ as can easily be seen.
\end{proof}

In regard to the previous proposition, we remark that the necessary
condition $F(t)=\wt{F}\iv(t)$ for the existence of an antisymmetric
factorization of $F\in G\cB\NN$ implies $F(1)^2=F(-1)^2=I$
by just putting $t=1$ or $t=-1$.
Hence for given $F$ (and thus given $F(1)$ and $F(-1)$),
the values of
\begin{equation}
\alpha+\beta,\qquad\gamma+\delta,\qquad
\alpha+\gamma,\qquad\beta+\delta.
\end{equation}
can immediately be determined without knowing the antisymmetric
factorization of $F$ or the corresponding characteristic pairs:

\subsection{Asymmetric factorizations within a Banach algebra}

In regard to the discussion at the end of Section \ref{s2} we are
going to show that each $A\in G\cB\NN$ can be factored in certain
``asymmetric'' ways.

As a first auxiliary step, we are going to specify the middle factors
$R(t)$, which appeared in (\ref{f2.2facts}) and (\ref{f2.2facts2}).
The following proposition shows the existence of
such factors, where the construction in the proof is completely explicit
(although not unique). Moreover, although we noted that the factors
$R(t)$ ought to be of a ``simple'' form, it turns out that the actually
important point is that they are related by means of the equation
$D(t)=R(t)W\wt{R}\iv(t)$ (in case of a factorization (\ref{f2.2facts}))
or the equation $D(t)=\wt{R}\iv(t)WR(t)$ (in case of a factorization 
(\ref{f2.2facts2})) to a factor $D(t)$ of the form (\ref{f3.D(t)}).

\begin{proposition}\label{p3.5}
Let $W\in\C\NN$ with $W^2=I$, and assume that $D(t)$ is given by
(\ref{f3.D(t)}) such that $D(1)\sim D(-1)\sim W$. Then there exists
a matrix function $R\in G\cB\NN$ such that $D(t)=R(t)W\wt{R}\iv(t)$.
\end{proposition}
\begin{proof}
We can assume that
\bqn
W &=& T\,\diag(I_{\sigma_+},-I_{\sigma_-})\,T\iv,
\eqn
where $\sigma_+$ and $\sigma_-$ are nonnegative integers with
$\sigma_++\sigma_-=N$ and $T\in G\C\NN$. With the numbers
$\alpha,\beta,\gamma,\delta$ defined as above in terms of the characteristic
pairs appearing in $D(t)$, it follows from Proposition \ref{p3.4}
(with $F(t)=D(t)$) that 
\begin{equation}
\sigma_+\;=\;\alpha+\beta\;=\;\alpha+\gamma\quad\mbox{ and }\quad
\sigma_-\;=\;\beta+\delta\;=\;\gamma+\delta.
\end{equation}
In particular, $\beta=\gamma$. Hence there exists a permutation matrix $\Pi_1$
such that
\bqn
D(t) &=& \Pi_1\,\diag(D_1(t),D_2(t),D_3(t))\,\Pi_1\iv,
\eqn
where
\bqn
D_1(t) &=&
\lefteqn{\under{\diag}{{\scriptstyle 0\le k\le \alpha}}\;(t^{\ka_k^{(1)}})}
\hspace{16ex}\mbox{ with $\ka_k^{(1)}\in\Z$ even,}
\\
D_2(t) &=&
\lefteqn{\under{\diag}{{\scriptstyle 0\le k\le \delta}}\;(-t^{\ka_k^{(2)}})}
\hspace{16ex}\mbox{ with $\ka_k^{(2)}\in\Z$ even,}
\\
D_3(t) &=&
\under{\diag}{{\scriptstyle 0\le k\le \beta}}\;
\left(\left(\ba{cc} t^{\ka_k^{(3)}} & 0 \\ 0 & -t^{\ka_k^{(4)}} \ea
\right)\right)\qquad
\mbox{ with $\ka_k^{(3)},\ka_k^{(4)}\in\Z$ odd.}
\eqn
Moreover, there exists another permutation matrix $\Pi_2$ such that
\bqn
W &=& T\Pi_2\,\diag(I_\alpha,-I_\delta,X_\beta)\;\Pi_2\iv T\iv,\\[1ex]
X_\beta &=& \under{\diag}{{\scriptstyle 1\le k\le\beta}}\;
\left(\left(\ba{cc} 1 & 0 \\ 0 & -1 \ea
\right)\right).
\eqn
We define 
\begin{equation}
R_1(t) \;=\;
\under{\diag}{{\scriptstyle 0\le k\le \alpha}}\;(t^{\ka_k^{(1)}/2}),
\qquad\qquad\qquad
R_2(t) \;=\;
\under{\diag}{{\scriptstyle 0\le k\le \delta}}\;(t^{\ka_k^{(2)}/2}),
\end{equation}
\bqn
R_3(t) &=&
\under{\diag}{{\scriptstyle 0\le k\le \beta}}\;
\left(\frac{1}{2}\left(\ba{cc}
t^{(\ka_k^{(3)}-\ka_k^{(4)})/2}+t^{(\ka_k^{(3)}+\ka_k^{(4)})/2} &
t^{(\ka_k^{(3)}-\ka_k^{(4)})/2}-t^{(\ka_k^{(3)}+\ka_k^{(4)})/2} \\
1-t^{\ka_k^{(4)}} & 1+ t^{\ka_k^{(4)}} \ea
\right)\right).\quad
\eqn
It can be verified straightforwardly that
\begin{equation}
D_1(t) \;=\; R_1(t)\wt{R}_1\iv(t),\qquad
D_2(t) \;=\; -R_2(t)\wt{R}_2\iv(t),\qquad
D_3(t) \;=\; R_3(t)X_\beta\wt{R}_3\iv(t).
\end{equation}
Hence, on defining $R$ by
\bqn
R(t) &=& \Pi_1\;\diag(R_1(t),R_2(t),R_3(t))\;\Pi_2\iv T\iv,
\eqn
it follows that $D(t)=R(t)W\wt{R}\iv(t)$.
\end{proof}

It is obvious that the previous proposition remains true if
one replaces the expression $D(t)=R(t)W\wt{R}\iv(t)$ with
$D(t)=\wt{R}\iv(t)WR(t)$. In fact, one just has to replace
$R(t)$ with $\wt{R}\iv(t)$, which is possible due to the assumption (c)
on the Banach algebra $\cB$.

Besides $\cB_+\NN$ and $\cB_-\NN$, we need another subalgebra of
$\cB\NN$. Given $W\in\C\NN$ with $W^2=I$, we define
\bqn\label{f2.BW}
\cB_W\NN &=& \cB\NN\cap(L^\iy)_W\NN.
\eqn
It is easy to see that $\cB_W\NN$ is an inverse closed Banach subalgebra
of $\cB\NN$.

The following theorem establishes the existence of two kinds of
``asymmetric'' factorizations within the Banach algebra $\cB$
for a given function $A\in G\cB\NN$.

\begin{theorem}\label{t3.6}
Let $W\in\C\NN$ with $W^2=I$, and assume that $A\in G\cB\NN$. Then
\begin{itemize}
\item[(a)]
there exists a factorization of $A(t)$ in the form
\bqn\label{f3.z2}
A(t) &=& A_-(t)R(t)A_0(t),
\eqn
where $A_-\in G\cB_-\NN$, $A_0\in G\cB_W\NN$, and $R\in G\cB\NN$ such that
$D(t)=R(t)W\wt{R}\iv(t)$ is of the form (\ref{f3.D(t)}).
Moreover,
\bqn\label{f3.z1}
F(t) &=& A_-(t)D(t)\wt{A}_-\iv(t)
\eqn
represents an antisymmetric factorization of the function
$F(t)=A(t)W\wt{A}\iv(t)$.
\item[(b)]
there exists a factorization of $A(t)$ in the form
\bqn\label{f3.z4}
A(t) &=& A_0(t)R(t)A_+(t),
\eqn
where $A_0\in G\cB_W\NN$, $A_-\in G\cB_-\NN$, and $R\in G\cB$ such that
$D(t)=\wt{R}\iv(t)WR(t)$ is of the form (\ref{f3.D(t)}).
Moreover,
\bqn\label{f3.z3}
G(t) &=& \wt{A}_+\iv(t)D(t)A_+(t)
\eqn
represents an antisymmetric factorization of the function
$G(t)=\wt{A}\iv(t)WA(t)$.
\end{itemize}
\end{theorem}
\begin{proof}
From the definition of $F$ it follows that
$F\in G\cB\NN$ and $F(t)=\wt{F}\iv(t)$. By Theorem \ref{t3.2}
there exists an antisymmetric factorizations (\ref{f3.z1}) with
$A_-\in G\cB_-\NN$ and $D(t)$ of the form (\ref{f3.D(t)}).

From the definition of $F$ it follows furthermore that $F(1)\sim F(-1)\sim W$,
which in turn implies $D(1)\sim D(-1)\sim W$.
Using Proposition \ref{p3.5} we obtain the existence of a function
$R\in G\cB\NN$ for which $D(t)=R(t)W\wt{R}\iv(t)$. Now we define
\bqn
A_0(t) &=& R\iv(t)A_-\iv(t)A(t),\nn
\eqn
which implies immediately the validity of equation (\ref{f3.z2}). Moreover,
\bqn
W\wt{A}_0(t)W &=& W\wt{R}\iv(t)\wt{A}_-\iv(t)\wt{A}(t)W.\nn
\eqn
From
$$
A(t)W\wt{A}\iv(t)\;=\;F(t)\;=\;A_-(t)D(t)\wt{A}_-\iv(t)
\;=\; A_-(t)R(t)W\wt{R}\iv\wt{A}_-\iv(t)
$$
we obtain
$$
W\wt{A}_0(t)W \;=\; R\iv(t)A_-\iv(t)A(t)\;=\;A_0(t).
$$
Hence $A_0\in \cB_W\NN$. Because, obviously, $A_0\in G\cB\NN$,
we even have $A_0\in G\cB_W\NN$. This settles part (b).

The proof of part (b) is similar. By Theorem \ref{t3.2} (see also the
remark made afterwards) and the facts that $G\in G\cB\NN$ and
$G(t)=\wt{G}\iv(t)$,
there exists an antisymmetric factorization (\ref{f3.z3}). From
$G(1)\sim G(-1)\sim W$ we obtain $D(1)\sim D(-1)\sim W$.
We apply again Proposition \ref{p3.5}, but now with $R$ replaced by
$\wt{R}\iv$, in order to conclude the existence of a function $R\in G\cB\NN$
for which $D(t)=\wt{R}\iv(t)WR(t)$. Finally, we define
\bqn
A_0(t) &=& A(t)A_+\iv(t)R\iv(t),\nn
\eqn
apply the equation
$$
\wt{A}\iv(t)WA(t) \;=\; G(t)\;=\;
\wt{A}_+\iv(t)D(t)A_+(t)\;=\;\wt{A}_+\iv(t)\wt{R}\iv(t)WR(t)A_+(t)
$$
and obtain in this way that $W\wt{A}_0(t)W=A_0(t)$. Hence $A_0\in G\cB_W\NN$.
\end{proof}

The proof of the previous theorem reveals that the
asymmetric factorization (\ref{f3.z2}) of $A(t)$ can be constructed in an
explicit way from the antisymmetric factorization of $F(t)$. Moreover,
to each possible middle factor $D(t)$ (hence, to each possible collection of
characteristic pairs), one can assign a corresponding function $R(t)$
which may appear as the middle factor in the asymmetric factorization.
The fact that this assignment may be carried out in different ways
does not affect the following considerations.

Similar statements hold, of course, also for the asymmetric factorization
(\ref{f3.z4}) of $A(t)$, which is connected with the antisymmetric
factorization of $G(t)$.

\section{Further properties of some classes of Toeplitz + Hankel operators}
\label{s4}

In this section we continue and conclude the study of the Toeplitz + Hankel operators
$\cM_W(A)$ and $\cN_W(A)$ with $A\in\cB\NN$. Notice first that it follows
from Proposition \ref{p2.2} and the inverse closedness of $\cB$ in
$C(\T)$ that these operators are Fredholm if and only if $A\in G\cB\NN$.

In the case $A\in G\cB\NN$ we will determine the dimension of the
kernel and cokernel of $\cM_W(A)$ and $\cN_W(A)$ in terms of
the characteristic pairs of an antisymmetric factorization of
a certain associated function. Formulas for the inverses (if they exist)
will also be presented.

For the following presentations, it is useful to introduce a
function $\Theta:\{-1,1\}\times\Z\to\Z$ which is defined by
\bqn
\Theta(\ro,\ka) &=& \left\{\ba{rr}
\ka/2 & \mbox{ if $\ka$ is even}\\
(\ka-\ro)/2 & \mbox{ if $\ka$ is odd}.\ea\right.
\eqn

For the interpretation of the following results, it is also helpful
to recall the notion of a pseudoinverse. Let $A$ be a linear
bounded operator acting on a Banach space $X$. A linear bounded
operator $A\kr$ acting also on $X$ is called a {\em pseudoinverse}
of $A$ if the relations
\begin{equation}
AA\kr A\;=\;A\quad\mbox{ and }\quad
A\kr AA\kr\;=\; A\kr
\end{equation}
hold. One can show that a pseudoinverse of $A$ exists if and only if
the image of $A$ is a complemented subspace in $X$, i.e., there exist
a closed subspace $X_0$ of $X$ such that $X={\rm im\,} A\oplus X_0$.
Hence each Fredholm operator possesses a pseudoinverse.
Pseudoinverses are in general in not unique. However, if $A$ is invertible,
then $A\kr$ is uniquely determined and coincides with $A\iv$.

\begin{lemma}\label{l4.1}
Let $D(t)$ be a matrix function of the form (\ref{f3.D(t)}) with the
characteristic pairs (\ref{f3.pairs1}). Then $H(D)^*=H(D)$ and
\bqn
\dim\ker(I+H(D)) &=&
\sum_{\ka_k>0}\Theta(\ro_k,\ka_k).
\eqn
\end{lemma}
\begin{proof}
First observe that $D^*(t)=\wt{D}(t)$. Consequently,
$$
H(D)^*\;=\;(PM(D)JP)^*\;=\;PJM(D^*)P\;=\;PM(\wt{D}^*)JP\;=\;H(D).
$$
Moreover, because
$H(D)=\diag(H(\ro_1t^{\ka_1}),H(\ro_2t^{\ka_2}),\dots,H(\ro_Nt^{\ka_N}))$
is a diagonal operator it suffices to determine
$\dim\ker(I+H(\ro_kt^{\ka_k}))$ and to
take the sum.  If $\ka_k\le0$, then $H(\ro_kt^{\ka_k})=0$.  Hence the
corresponding dimension is zero.  If $\ka_k>0$, then the matrix
representation of $H(\ro_kt^{\ka_k})$ has entries $\ro_k$ only on the
$\ka_k$-th diagonal, which connects the entries $(1,\ka_k)$ and
$(\ka_k,1)$, and has zero entries elsewhere.  From this it is easy to see
that the dimension equals $\Theta(\ro_k,\ka_k)$. 
\end{proof}

\begin{lemma}\label{l4.2}
Assume that $D(t)$ is of the form (\ref{f3.D(t)}). Then
\begin{equation}\label{f4.HD}
T(\wt{D})H(D)\;=\;H(D)T(D)\;=\;0\quad\mbox{ and }\quad
H(D)^3\;=\;H(D).
\end{equation}
Moreover, if we introduce
\bqn\label{f4.BBkr}
&&\lefteqn{B\;=\;I+H(D),}\hspace{22ex}
B\kr\;=\;I-H(D)^2+\textstyle\frac{1}{4}(H(D)^2+H(D)),
\eqn
then $B\kr B B\kr=B\kr$ and $BB\kr B=B$.
\end{lemma}
\begin{proof}
By considering the scalar case, $D(t)=\ro_kt^{\ka_k}$ and distinguishing
$\ka_k>\le0$ and $\ka_k>0$, it can be seen straightforwardly
that $T(\wt{D})H(D)=H(D)T(D)=0$.
Moreover, using (\ref{f.Tab}) and the fact that $D(t)=\wt{D}\iv(t)$
it follows that
$$
H(D)^3\;=\;H(D)H(D)H(\wt{D}\iv)\;=\;H(D)(I-T(D)T(D\iv))\;=\;H(D).
$$
In order to prove that $B\kr B B\kr=B\kr$ and $BB\kr B=B$, we introduce
$p=I-H(D)^2$ and $q=(H(D)+H(D)^2)/2$. By just using the identity
$H(D)^3=H(D)$, one can verify that $p^2=p$, $q^2=q$, $pq=qp=0$.
Because $B=p+2q$ and $B\kr=p+q/2$, the desired relations follow
immediately.
\end{proof}

\begin{lemma}\label{l4.3}
Let $W\in\C\NN$ with $W^2=I$, and assume that $R\in G\cB\NN$ is
given such that $D(t)=R(t)W\wt{R}\iv(t)$ is of the form (\ref{f3.D(t)}).
Introduce the operators $B$ and $B\kr$ by (\ref{f4.BBkr}), and
\bqn
&&\lefteqn{A_1\;=\;\cM_W(R),}   \hspace{25ex}A_2\;=\;\cN_W(R\iv),\\
&&\lefteqn{A_1\kr\;=\; A_2B\kr,}\hspace{25ex}A_2\kr\;=\;B\kr A_1.
\eqn
Then $B=A_1A_2$, $A_1\kr=\cN_W(R\iv)(I-\frac{1}{2}H(D))$,
$A_2\kr=(I-\frac{1}{2}H(D))\cM_W(R)$, and
\begin{equation}
A_1\kr A_1A_1\kr\;=\; A_1\kr,\quad
A_1A_1\kr A_1\;=\; A_1,\quad
A_2\kr A_2A_2\kr\;=\; A_2\kr,\quad
A_2A_2\kr A_2\;=\; A_2.
\end{equation}
\end{lemma}
\begin{proof}
Using formula (\ref{f2.MWNW}), it follows that
$$
\cM_W(R)\cN_W(R\iv)\;=\;T(RR\iv)+H(RW\wt{R}\iv)
\;=\;I+H(D).
$$
Hence $A_1A_2=B$.
By using this, the relations of $A_1\kr=A_2B\kr$ and $A_2\kr=B\kr A_1$,
and formula $B\kr BB\kr=B\kr$ from Lemma \ref{l4.2}, we obtain
\bqn
&&A_1\kr A_1A_1\kr \;=\; A_2B\kr A_1A_2B\kr \;=\; A_2B\kr BB\kr
\;=\;A_2 B\kr\;=\;A_1\kr,\nn \\
&&A_2\kr A_2A_2\kr \;=\; B\kr A_1A_2B\kr A_1 \;=\; B\kr BB\kr A_1
\;=\;B\kr A_1\;=\;A_2\kr.\nn
\eqn

Next, as an auxiliary step, we are going to establish the identities
\bqn\label{f4.AHH}
H(D)A_1 \;=\; H(D)^2A_1\quad\mbox{ and }\quad
A_2H(D) \;=\; A_2H(D)^2.
\eqn
Indeed, using that $R=D\wt{R}W$ and formulas (\ref{f.Tab}) and (\ref{f.Hab}),
it follows that
\bqn
A_1 &=& T(D\wt{R}W)+H(D\wt{R})\nn\\
&=& T(D)T(\wt{R}W)+H(D)H(RW)+T(D)H(\wt{R})+H(D)T(R)\nn\\
&=& H(D)A_1+T(D)(T(\wt{R}W)+H(\wt{R})).\nn
\eqn
Multiplying from the left with $H(D)$ and observing that $H(D)T(D)=0$ by
Lemma \ref{l4.2}, we obtain the first identity in (\ref{f4.AHH}).
Similarly, by using $R\iv=W\wt{R}\iv\wt{D}$, it follows that
\bqn
A_2 &=& T(W\wt{R}\iv\wt{D})+H(R\iv D)\nn\\
&=& T(W\wt{R}\iv)T(\wt{D})+H(W\wt{R}\iv)H(D)+T(R\iv)H(D)+H(R\iv)T(\wt{D})\nn\\
&=& A_2H(D)+(T(W\wt{R}\iv)+H(R\iv))T(\wt{D}).\nn
\eqn
Multiplying from the right with $H(D)$ by observing that $T(\wt{D})H(D)=0$,
we arrive at the second identity in (\ref{f4.AHH}).

Using this and the definition of $B\kr$, it follows that
$A_1\kr=A_2(I-\frac{1}{2}H(D))$ and $A_2\kr=(I-\frac{1}{2}H(D))A_1$.
Hence we obtain the desired expressions for $A_1\kr$ and $A_2\kr$.

Moreover, using the notation $p$ and $q$ introduced in the proof of
Lemma \ref{l4.2}, it is easy to see that $BB\kr=B\kr B=p+q$. Hence
\begin{equation}
BB\kr\;=\;B\kr B\;=\; I+\textstyle\frac{1}{2}(H(D)-H(D)^2).
\end{equation}
Combining this with (\ref{f4.AHH}) it follows that
$$
BB\kr A_1\;=\;A_1\quad\mbox{ and }\quad
A_2B\kr B\;=\;A_2.
$$
Now we are able to derive the remaining identities:
\bqn
&&A_1A_1\kr A_1\;=\; A_1A_2B\kr A_1\;=\;BB\kr A_1\;=\;A_1,\nn\\
&&A_2A_2\kr A_2\;=\; A_2B\kr A_1A_2\;=\;A_2B\kr B\;=\;A_2.\nn
\eqn
This completes the proof.
\end{proof}

\begin{lemma}\label{l4.4}
Let $A_1$, $A_2$ and $B$ as before. Then
\begin{equation}
\dim\ker A_1^*\;=\;\dim\ker A_2\;=\;\dim\ker B
\end{equation}
\end{lemma}
\begin{proof}
Since $H(D)^*=H(D)$ by Lemma \ref{l4.1}, it follows that $B^*=B$.
The relation $B=A_1A_2$ stated in Lemma \ref{l4.3}
implies that $\ker A_2\subseteq\ker B$ and
$\ker A_1^*\subseteq\ker B^*$.
Moreover, because $A_2=A_2A_2\kr A_2=A_2B\kr A_1A_2=A_1\kr B$, we obtain
$\ker B\subseteq\ker A_2$. Similarly, since $A_1=A_1A_1\kr A_1=
A_1A_2B\kr A_1=BA_2\kr$, we arrive at $\ker B^*\subseteq \ker A_1^*$.
\end{proof}

Now we are able to establish formulas for the dimension of the kernel
and cokernel of the operators $\cM_W(R)$ and $\cN_W(R)$, where 
$R(t)$ represent appropriate middle factors that are expected to appear
in the asymmetric factorization. Notice the slightly modified
notation in the following proposition, i.e., we are considering
$\cN_W(R\iv)$ instead of $\cN_W(R)$. The important point is, however, that
$R(t)$ is related to a matrix function $D(t)$ of the form (\ref{f3.D(t)}).

\begin{proposition}\label{p4.5}
Let $W\in\C\NN$ with $W^2=I$, and assume that $R\in G\cB\NN$
is given such that $D(t)=R(t)W\wt{R}\iv(t)$ is of the form
(\ref{f3.D(t)}) with characteristic pairs (\ref{f3.pairs1}). Then 
\bqn
\dim\ker\cM_W(R) &=& \lefteqn{-\sum_{\ka_k<0}\Theta(\ro_k,\ka_k),}
\hspace{37ex}
\makebox[0ex][r]{$\dim\ker\cM_W(R)^*$} \;=\;
\sum_{\ka_k>0}\Theta(\ro_k,\ka_k),
\\
\dim\ker\cN_W(R\iv) &=& \lefteqn{\sum_{\ka_k>0}\Theta(\ro_k,\ka_k),}
\hspace{37ex}
\makebox[0ex][r]{$\dim\ker\cN_W(R\iv)^*$} \;=\;
-\sum_{\ka_k<0}\Theta(\ro_k,\ka_k).\qquad
\eqn
\end{proposition}
\begin{proof}
The formulas for $\dim\ker\cM_W(R)^*$ and $\dim\ker\cN_W(R\iv)$ follow
immediately from Lemma \ref{l4.4} in connection with Lemma \ref{l4.1}.
Moreover, by the index formula stated in Proposition \ref{p2.2}
it can be seen that
\bqn
\ind\cM_W(R)&=&\ind T(R) \;=\; -\wind\det R,\nn\\
\ind\cN_W(R\iv)&=&\ind T(R\iv) \;=\; -\wind\det R\iv\;=\;\wind\det R.\nn
\eqn
Because $D(t)=R(t)W\wt{R}\iv(t)$, we have
$$
\wind \det D \;=\; \wind\det R+\wind\det \wt{R}\iv\;=\;2\,\wind\det R.
$$
On the other hand,
\bqn
\wind\det D &=& \sum_{k=1}^N \ka_k.\nn
\eqn
Combining all this, it follows that
\begin{equation}
\ind\cM_W(R) \;=\; -\frac{1}{2}\,\sum_{k=1}^N\ka_k,\qquad
\ind\cN_W(R\iv) \;=\;  \frac{1}{2}\,\sum_{k=1}^N\ka_k.
\end{equation}
Since $D(1)\sim D(-1)\sim W$, we obtain from Proposition \ref{p3.4},
in particular, that $\beta=\gamma$. We conclude from the definition
of the function $\Theta$ that
\bqn
\sum_{k=1}^N \Theta(\ro_k,\ka_k) &=& 
\frac{1}{2}\sum_{k=1}^N\ka_k-\frac{1}{2}\sum_{\ka_k{\rm\;odd}}\ro_k
\nn\\
&=& \frac{1}{2}\sum_{k=1}^N\ka_k+\frac{\gamma-\beta}{2}
\;=\;\frac{1}{2}\sum_{k=1}^N\ka_k.
\eqn
Using the formulas
\bqn
\ind\cM_W(R) &=& \dim\ker\cM_W(R)-\dim\ker\cM_W(R)^*,\nn\\
\ind\cN_W(R\iv) &=& \dim\ker\cN_W(R\iv)-\dim\ker\cN_W(R\iv)^*,\nn
\eqn
it is easy to derive the remaining two formulas.
\end{proof}

The following result determines the dimensions of the kernel and
cokernel of the operators $\cM_W(A)$ and $\cN_W(A)$ for $A\in G\cB\NN$
in terms of the characteristic pairs of an associated antisymmetric
factorization problem. Note that the existence of this antisymmetric
factorization is ensured by Theorem \ref{t3.6}.

Moreover, we give expressions for the pseudoinverses of the above
operator, which are the inverses in case of invertibility.  It should
also be observed that in the formulation of the following theorem we
need not make reference to the asymmetric factorizations, although
they are, of course, used in the proof.

\begin{theorem}\label{t4.6}
Let $W\in\C\NN$ with $W^2=I$, and let $A\in G\cB\NN$.
\begin{itemize}
\item[(a)]
Assume that an antisymmetric factorization of the
function $F(t)=A(t)W\wt{A}\iv(t)$ is given by
$F(t)=A_-(t)D(t)\wt{A}_-\iv(t)$, where $A_-\in G\cB_-\NN$ and
$D(t)$ is of the form (\ref{f3.D(t)}) with the
characteristic pairs (\ref{f3.pairs1}). Then
\bqn\label{f4.ker1}
\dim\ker\cM_W(A)   &=& -\sum_{\ka_k<0}\Theta(\ro_k,\ka_k),\\
\dim\ker\cM_W(A)^* &=& \sum_{\ka_k>0}\Theta(\ro_k,\ka_k).\label{f4.ker2}
\eqn
Moreover, a pseudoinverse of $\cM_W(A)$ is given by
\begin{equation}\label{f4.pseu1}
\cN_W(A\iv A_-)(I-\textstyle\frac{1}{2}H(D))T(A_-\iv).
\end{equation}
\item[(b)]
Assume that an antisymmetric factorization of the
function $G(t)=\wt{A}\iv(t)WA(t)$ is given by
$G(t)=\wt{A}_+\iv(t)D(t)A_+(t)$, where $A_+\in G\cB_+\NN$ and
$D(t)$ is of the form (\ref{f3.D(t)}) with the
characteristic pairs (\ref{f3.pairs1}). Then
\bqn\label{f4.ker3}
\dim\ker\cN_W(R)   &=& -\sum_{\ka_k<0}\Theta(-\ro_k,\ka_k),\\
\dim\ker\cN_W(R)^* &=& \sum_{\ka_k>0}\Theta(-\ro_k,\ka_k).\label{f4.ker4}
\eqn
Moreover, a pseudoinverse of $\cN_W(A)$ is given by
\begin{equation}\label{f4.pseu2}
T(A_+\iv)(I-\textstyle\frac{1}{2}H(D\iv))\cM_W(A_+A\iv).
\end{equation}
\end{itemize}
\end{theorem}
\begin{proof}
Let us first consider case (a). By Theorem \ref{t3.6} we can assume that
we are given an asymmetric factorization $A(t)=A_-(t)R(t)A_0(t)$
with the conditions on the factors stated there. In addition, we
are given an antisymmetric factorization $F(t)=A_-(t)D(t)\wt{A}_-\iv(t)$
with $D(t)=R(t)W\wt{R}\iv(t)$ of the function $F(t)=A(t)W\wt{A}\iv(t)$.
From (\ref{f2.Mmult}) and (\ref{f2.M=T}) it follows that
\bqn
\cM_W(A) &=& \cM_W(A_-)\cM_W(R)\cM_W(A_0),\nn
\eqn
where both $\cM_W(A_-)=T(A_-)$ and $\cM_W(A_0)$ are invertible.
There inverses are equal to $T(A_-\iv)$ and $\cM_W(A_0\iv)$, respectively.
Hence the dimension of the kernel and cokernel of $\cM_W(A)$ is equal to
that of $\cM_W(R)$, which, in turn, has been given in
Proposition \ref{p4.5}.
	
Next we need to take into account the following fact, which can be proved
straightforwardly: if an operator $S_1$ has a pseudoinverse $S_1\kr$
and $S_2=US_2V$ where $U$ and $V$ are invertible operators, then
a pseudoinverse of $S_2$ is given by $S_2\kr=V\iv S_1\kr U\iv$.

It follows from Lemma \ref{l4.3} that a pseudoinverse of $\cM_W(R)=A_1$
is given by $\cN_W(R\iv)(I-\frac{1}{2}H(D))=A_1\kr$. Consequently, a
pseudoinverse of $\cM_W(A)$ is given by
$$\textstyle
\cM_W(A_0\iv)\cN_W(R\iv)(I-\frac{1}{2}H(D))T(A_-\iv).
$$
Now we use formula (\ref{f2.M=N}) and (\ref{f2.Nmult}) in order to conclude
that
$$\cM_W(A_0\iv)\cN_W(R\iv)\;=\; \cN_W(A_0\iv)\cN_W(R\iv)\;=\;
\cN_W(A_0\iv R\iv).$$
Remark that $A_0\iv R\iv=A\iv A_-$ (because this is just the
equation $A=A_-RA_0$). Combining these last facts, we arrive at the
desired expression for the pseudoinverse.

Case (b) can be treated in the same way, but we give the complete proof because
the notation differs sometimes here in comparison with previous results.
First of all, we may assume that we are given an asymmetric factorization
$A(t)=A_0(t)R(t)A_+(t)$ as has been stated in Theorem \ref{t3.6}.
Moreover, we are given an antisymmetric factorization
$G(t)=\wt{A}_+\iv(t)D(t)A_+(t)$ with $D(t)=\wt{R}\iv(t)WR(t)$ of the
function $G(t)=\wt{A}\iv(t)WA(t)$. From (\ref{f2.Nmult}) and (\ref{f2.N=T})
it follows that
\bqn
\cN_W(A) &=& \cN_W(A_0)\cN_W(R)\cN_W(A_+),\nn
\eqn
where $\cN_W(A_0)$ and $\cN_W(A_+)=T(A_+)$ are invertible.
The inverses are $\cN_W(A_0\iv)$ and $T(A_+\iv)$, respectively.
The above relation $D(t)=\wt{R}\iv(t)WR(t)$ can be rewritten
as $D\iv(t)=R\iv(t)W\wt{R}(t)$.
We now have to apply Proposition \ref{p4.5} and Lemma \ref{l4.3}
with $R(t)$ replaced with $R\iv(t)$ and $D(t)$ replaced with $D\iv(t)$.
Correspondingly, the characteristic pairs $(\ro_k,\ka_k)$ have to be
replaced with $(\ro_k,-\ka_k)$. We arrive at the formulas
$$
\dim\ker\cN_W(R) \;=\; \sum_{\ka_k<0}\Theta(\ro_k,-\ka_k),\qquad
\dim\ker\cN_W(R)^* \;=\;
-\sum_{\ka_k>0}\Theta(\ro_k,-\ka_k).
$$
It remains to note that $\Theta(\ro_k,-\ka_k)=-\Theta(-\ro_k,\ka_k)$
in order to conclude the desired formulas for the dimension of the
kernel and cokernel of $\cN_W(A)$.

Also in regard to the pseudoinverse we have to apply Lemma \ref{l4.3}, but
with $R(t)$ replaced with $R\iv(t)$ and $D(t)$ replaced with $D\iv(t)$.
It follows that the pseudoinverse of $\cN_W(R)=A_2$ is given by
$\textstyle (I-\frac{1}{2}H(D\iv))\cM_W(R\iv)=A_2\kr$.
As before, we obtain that a pseudoinverse of
$\cN_W(A)$ is given by
$$
T(A_+\iv)\textstyle (I-\frac{1}{2}H(D\iv))\cM_W(R\iv)\cN_W(A_0\iv).
$$
Using formulas (\ref{f2.M=N}) and (\ref{f2.Mmult}) we derive
$$
\cM_W(R\iv)\cN_W(A_0\iv)\;=\;\cM_W(R\iv)\cM_W(A_0\iv)\;=\;\cM_W(R\iv A_0\iv).
$$
The desired pseudoinverse of $\cM_W$ is now obtained by piecing together these
last facts in connection with $R\iv A_0\iv=A_+A\iv$, which is just
the factorization $A=A_0RA_+$ rewritten.
\end{proof}

At the end of this section we consider some simple consequences of the
previous theorem. In particular, we state the necessary and sufficient
conditions for the invertibility of the operators $\cM_W(A)$ and $\cN_W(A)$.

\begin{corollary}\label{c4.7}
Let $W\in\C\NN$ with $W^2=I$, and assume $A\in G\cB\NN$. 
\begin{itemize}
\item[(a)] 
The operator $\cM_W(A)$ is invertible if and only if the function
$F(t)=A(t)W\wt{A}\iv(t)$ admits an antisymmetric factorization with
characteristic pairs $(\ro_k,\ka_k)$
which are all contained in the set
\begin{equation}\label{f4.set1}
\Big\{\;(-1,-1),\,(-1,0),\,(1,0),\,(1,1)\;\Big\}.
\end{equation}
\item[(b)] 
The operator $\cN_W(A)$ is invertible if and only if the function
$G(t)=\wt{A}\iv(t)WA(t)$ admits an antisymmetric factorization with
characteristic pairs $(\ro_k,\ka_k)$
which are all contained in the set
\begin{equation}\label{f4.set2}
\Big\{\;(1,-1),\,(-1,0),\,(1,0),\,(-1,1)\;\Big\}.
\end{equation}
\end{itemize}
\end{corollary}
\begin{proof}
The operators are invertible if and only if the sums in (\ref{f4.ker1})
and (\ref{f4.ker2}), or, (\ref{f4.ker3}) and (\ref{f4.ker4}), respectively
are zero. Notice that the different terms appearing there are all nonnegative
integers. Hence they must be equal to zero. It remains to remark that
$\Theta(\ro,\ka)=0$ if and only if $\ka=0$ or $\ka=\ro=1$ or $\ka=\ro=-1$.
\end{proof}

The previous result takes a much simpler form in the two special cases
where $W=I$ or $W=-I$. In fact, we can apply Proposition \ref{p3.4}
and recall the definition of the numbers $\alpha,\beta,\gamma,\delta$.

In the case where $W=I$, we have $F(1)=F(-1)=I$ and $G(1)=G(-1)=I$.
Hence Proposition \ref{p3.4} implies that $\alpha=N$ and
$\beta=\gamma=\delta=0$. Hence among the pairs given in (\ref{f4.set1})
or (\ref{f4.set2}) only the pair $(1,0)$ can occur. The result is that
the operator $\cM_I(A)$ ($\cN_I(A)$, resp.) is invertible if and only if
the function $F(t)=A(t)\wt{A}\iv(t)$ ($G(t)=\wt{A}\iv(t)A(t)$, resp.)
admits an antisymmetric factorization with all characteristic pairs
equal to $(1,0)$.

In the case where $W=-I$, we obtain in a similar way the result that
the operator $\cM_{-I}(A)$ ($\cN_{-I}(A)$, resp.) is invertible if and only if
the function $F(t)=-A(t)\wt{A}\iv(t)$ ($G(t)=-\wt{A}\iv(t)A(t)$, resp.)
admits an antisymmetric factorization with all characteristic pairs
equal to $(-1,0)$.

\section{Singular integral operators with flip}

In this section, we study the properties of singular integral
operators with flip. In particular, we obtain results for the dimension
of the kernel and cokernel in the case of Fredholmness under the assumption
that the generating functions belongs to the Banach algebra $\cB\NN$.

In what follows, when we are given the matrix functions
$a,b,c,d\in L^\iy(\T)\NN$, we associate a matrix function 
of twice the matrix size,
\bqn\label{f5.1}
A &=& \left(\ba{cc} a & b \\ c & d \ea\right)
\in L^\iy(\T)\NNN.
\eqn
Moreover, we let $W$ stand for the following constant matrix of size
$2N\times 2N$,
\bqn\label{f5.W}
W &=& \left(\ba{cc} 0 & I_N \\ I_N & 0 \ea\right).
\eqn
Finally, to a matrix function $A\in  L^\iy(\T)\NNN$ given as above, we
associate a matrix function $\wh{A}\in  L^\iy(\T)\NNN$ defined by
\bqn
\wh{A}(t) &=& W\wt{A}(t)W.
\eqn

Next we introduce, for given $A\in  L^\iy(\T)\NNN$, the following operators
\bqn
\label{f5.cT}
\cT(A) &=& \Big(P,\,JP\Big)M(A)\left(\ba{c} P\\ PJ\ea\right),\\
\label{f5.cH}
\cH(A) &=& \Big(P,\,JP\Big)M(A)\left(\ba{c} Q\\ QJ\ea\right).
\eqn
Using the basic relations for the operators $P$, $Q$, $J$, and $M(A)$,
it is easy to see that then
\bqn
\cT(\wh{A}) &=& \Big(Q,\,JQ\Big)M(A)\left(\ba{c} Q\\ QJ\ea\right),\\
\cH(\wh{A}) &=& \Big(Q,\,JQ\Big)M(A)\left(\ba{c} P\\ PJ\ea\right).
\eqn
Moreover, given $A,B\in L^\iy(\T)\NNN$, the following relations hold:
\bqn
\label{f5.cTab}
\cT(AB) &=& \cT(A)\cT(B)+\cH(A)\cH(\wh{B}),\\
\label{f5.cHab}
\cH(AB) &=& \cT(A)\cH(B)+\cH(A)\cT(\wh{B}).
\eqn
In fact, they are essentially a consequence of the identity
$$
\left(\ba{c} P\\ PJ\ea\right)\Big(P,\,JP\Big)+
\left(\ba{c} Q\\ QJ\ea\right)\Big(Q,\,JQ\Big)
\;=\; \left(\ba{cc} I&0\\ 0&I\ea\right)
$$
The formal resemblance to the formulas (\ref{f.Tab}) and (\ref{f.Hab})
is obvious.

Finally, we define the operators
\bqn
\label{f5.Phi}
\Phi(A) &=& \cT(A)+\cH(A) \;=\;
\Big(P,\,JP\Big)M(A)\left(\ba{c} I\\ J\ea\right),\\
\label{f5.Psi}
\Psi(A) &=& \cT(A)+\cH(\wh{A}) \;=\;
\Big(I,\,J\Big)M(A)\left(\ba{c} P\\ PJ\ea\right).
\eqn
These operators are the {\em singular integral operators with flip}
which we intend to study in this section. Indeed, if
$A$ is given by (\ref{f5.1}), then
\bqn
\label{f5.124}
\Phi(A) &=& 
PM(a)+PJM(\tilde{b})+QJM(c)+QM(\tilde{d}),\\
\label{f5.125}
\Psi(A) &=&
M(a)P+M(b)JQ+M(\tilde{c})JP+M(\tilde{d})Q.
\eqn

Using the formulas (\ref{f5.cTab}) and (\ref{f5.cHab}), formulas
analogous to (\ref{f.Mab}) and (\ref{f.Nab}) can be derived:
\bqn
\Phi(AB) &=& \Phi(A)\Phi(B)+\cH(A)\Phi(\wh{B}-B),\\
\Psi(AB) &=& \Psi(A)\Psi(B)+\Psi(\wh{A}-A)\cH(\wh{B}).
\eqn
Indeed, 
\bqn
\Phi(AB) &=& \cT(AB)+\cH(AB)\nn\\
&=& \cT(A)\cT(B)+\cH(A)\cH(\wh{B})+\cT(A)\cH(B)+\cH(A)\cT(\wh{B})\nn\\
&=& \cT(A)\Phi(B)+\cH(A)\Phi(\wh{B})\nn\\
&=& \Phi(A)\Phi(B)+\cH(A)\Phi(\wh{B}-B).\nn
\eqn
Moreover,
\bqn
\Psi(AB) &=& \cT(AB)+\cH(\wh{A}\wh{B})\nn\\
&=& \cT(A)\cT(B)+\cH(A)\cH(\wh{B})+\cT(\wh{A})\cH(\wh{B})+
\cH(\wh{A})\cT(B)\nn\\
&=& \Psi(A)\cT(B)+\Psi(\wh{A})\cH(\wh{B})\nn\\
&=& \Psi(A)\Psi(B)+\Psi(\wh{A}-A)\cH(\wh{B}).\nn
\eqn

The corresponding ``simplifications'', where multiplicativity holds,
read as follows:
\bqn
\Phi(AB) &=&\lefteqn{\Phi(A)\Phi(B)}\hspace{15ex}
\mbox{if $A\in (\ovl{H^\iy})\NNN$ or $B\in (L^\iy)\NNN_W$};
\label{f5.Phimult}\\
\Psi(AB) &=&\lefteqn{\Psi(A)\Psi(B)}\hspace{15ex}
\mbox{if $A\in (L^\iy)\NNN_W$ or $B\in (H^\iy)\NNN$}.
\label{f5.Psimult}
\eqn
Here $W$ is given by (\ref{f5.W}). Notice that the Banach algebra
$(L^\iy)\NNN_W$ is equal to
\bqn
\Big\{\;A\in L^\iy(\T)\NNN\;:\; \wh{A}=A\;\Big\}.
\eqn

Also the counterpart to formula (\ref{f2.MWNW}) can be established:
\bqn\label{f5.PhiPsi}
\Phi(A)\Psi(B) &=& \cT(AB)+\cH(A\wh{B}).
\eqn
Indeed, using (\ref{f5.cTab}) and (\ref{f5.cHab}) it follows that
\bqn
\Phi(A)\Psi(B) &=&
\Big(\cT(A)+\cH(A)\Big)\Big(\cT(B)+\cH(\wh{B})\Big)\nn\\
&=& \cT(A)\cT(B)+\cH(A)\cH(\wh{B})+\cH(A)\cT(B)+\cT(A)\cH(\wh{B})\nn\\
&=& \cT(AB)+\cH(A\wh{B}).\nn
\eqn

The analogy of these formulas in comparison with previous formulas
finds its crystal explanation in the following result.

\begin{proposition}\label{p5.1}
The mapping $\Xi$ defined by
\bqn
\Xi &:& \cL(L^2)\NN\to\cL(H^2)\NNN,\;\;
X\mapsto\left(\ba{c} P\\ PJ\ea\right)X\Big(P,\;JP\Big)
\eqn
represents a C*-algebra isomorphism between $\cL(L^2)\NN$ and
$\cL(H^2)\NNN$. In particular, the mapping $\Xi$ acts as follows:
\bqn
\label{f5.Xi1}
&&\lefteqn{\Xi:\cT(A)\mapsto T(A),}\hspace{30ex}
\Xi:\cH(A)\mapsto H(AW),\\
\label{f5.Xi2}
&&\lefteqn{\Xi:\Phi(A)\mapsto \cM_W(A),}\hspace{30ex}
\Xi:\Psi(A)\mapsto \cN_W(A),
\eqn
for $A\in L^\iy(\T)\NNN$, where $W$ is given by (\ref{f5.W}).
\end{proposition}
\begin{proof}
The first assertion follows from the fact that the linear operators
\bqn
\Big(P,\;JP\Big):(H^2)^{2N}\to(L^2)^N\quad\mbox{ and }\quad
\left(\ba{c} P\\ PJ\ea\right):(L^2)^N\to(H^2)^{2N}
\eqn
are Hilbert space isometries and are both the inverse and the
adjoint of each other. In fact,
\bqn
\Big(P,\;JP\Big)\left(\ba{c} P\\ PJ\ea\right) &=& I,\nn\\
\left(\ba{c} P\\ PJ\ea\right)\Big(P,\;JP\Big) &=&
\left(\ba{cc} P&0\\ 0&P\ea\right).\nn
\eqn
In order to prove (\ref{f5.Xi1}) and (\ref{f5.Xi2})
it suffices to recall the definitions (\ref{f5.cT}), (\ref{f5.cH}),
(\ref{f5.Phi}), (\ref{f5.Psi}), to use the last identity and
the relation
\bqn
\left(\ba{c} Q\\ QJ\ea\right) &=&
\left(\ba{cc} 0 & J \\  J & 0 \ea\right)
\left(\ba{c} P\\ PJ\ea\right)
\;\;=\;\;
W\left(\ba{cc} J&0\\ 0&J\ea\right)
\left(\ba{c} P\\ PJ\ea\right).\nn
\eqn
One has also to use the definition of $\cM_W(A)$ and $\cN_W(A)$
and the fact that $\wh{A}=W\wt{A}W$.
\end{proof}

The importance of the previous proposition is that it says that
the above singular integral operators with flip are unitarily
equivalent to Toeplitz + Hankel operators. These Toeplitz + Hankel
operators fall exactly in the classes which were studied in the previous
sections. Hence it is possible to reduce the study of the several properties
of singular integral operators with flip to the corresponding problems
for these Toeplitz + Hankel operators, which has already been done.

In regard to Proposition \ref{p2.2} and Proposition \ref{p2.4},
the following result is an immediate consequence.

\begin{proposition}\label{p5.2}
Let $A\in L^\iy(\T)\NNN$.
\begin{itemize}
\item[(a)]
If $\Phi(A)$ is Fredholm, then $A\in G(L^\iy(\T)\NNN)$.
\item[(b)]
If $\Psi(A)$ is Fredholm, then $A\in G(L^\iy(\T)\NNN)$.
\end{itemize}
Now assume that $A\in C(\T)\NNN$. Then
\begin{itemize}
\item[(c)]
$\Phi(A)$ is Fredholm if and only if $A\in G(C(\T)\NNN)$.
\item[(d)]
$\Psi(A)$ is Fredholm if and only if $A\in G(C(\T)\NNN)$.
\end{itemize}
Moreover, if this is true, then $\ind\Phi(A)=\ind\Psi(A)=-\wind\det A$.
\end{proposition}

Another result concerns the case where $A\in(L^\iy)_W\NNN$
with $W$ given by (\ref{f5.W}). As we will see shortly, this case is trivial.
If $A$ is given by (\ref{f5.1}) and $\wh{A}=A$, then
$d=\tilde{a}$ and $c=\tilde{b}$. In other words,
\bqn\label{f5.MMJ}
\Phi(A) \;\;=\;\; \Psi(A) \;\;=\;\;
M(a)+M(b)J,
\eqn
which is an operator composed of multiplication operators and 
the flip operator, but without the usual singular integral operator $S=P-Q$.
Compare in this connection the first equality in this formula
with the identity (\ref{f2.M=N}).

For completeness sake, we state the corresponding
invertibility and Fredholm criteria for operators (\ref{f5.MMJ}),
which follow from Corollary \ref{c2.3} by means
of Proposition \ref{p5.1}. It can be proved also by different,
more straightforward considerations.

\begin{corollary}\label{c5.3}
Let $A\in(L^\iy)_W\NNN$, where $W$ is given by (\ref{f5.W}). Then
the following is equivalent:
\begin{itemize}
\item[(i)]
$A\in G(L^\iy)_W\NNN$.
\item[(ii)]
$\Phi(A)=\Psi(A)$ is invertible.
\item[(iii)]
$\Phi(A)=\Psi(A)$ is Fredholm.
\end{itemize}
If this is fulfilled, then the inverse of $\Phi(A)=\Psi(A)$ is given by
$\Phi(A\iv)=\Psi(A\iv)$.
\end{corollary}

Now we turn to the case in which we are actually interestated in,
namely the  operators $\Phi(A)$ and $\Psi(A)$ with $A\in\cB\NNN$.
As before we assume that the Banach algebra $\cB$ possesses the properties
(a)--(d) stated at the beginning of Section \ref{s3}.

Because the Banach algebra $\cB$ is inverse closed in $C(\T)$,
Proposition \ref{p5.2}(cd) implies that, for given $A\in\cB\NNN$,
the operator $\Phi(A)$ ($\Psi(A)$, resp.) is a Fredholm operator if
and only if $A\in G\cB\NNN$. Similar as in Section \ref{s4}, and, of course,
referring to these results, we will determine the dimension of the kernel
and cokernel of $\Phi(A)$ and $\Psi(A)$ in terms of the characteristic pairs
of an antisymmetric factorization of a certain associated function.
Formulas for pseudoinverses (which are the inverses in the case
of invertibility) will also be presented.

\begin{theorem}\label{t5.4}
Let $W$ be given by (\ref{f5.W}), and let $A\in G\cB\NNN$.
\begin{itemize}
\item[(a)]
Assume that an antisymmetric factorization of the
function $F(t)=A(t)W\wt{A}\iv(t)$ is given by
$F(t)=A_-(t)D(t)\wt{A}_-\iv(t)$, where $A_-\in G\cB_-\NN$ and
$D(t)$ is of the form (\ref{f3.D(t)}) with the
characteristic pairs (\ref{f3.pairs1}). Then
\bqn\label{f5.ker1}
\dim\ker\Phi(A)   &=& -\sum_{\ka_k<0}\Theta(\ro_k,\ka_k),\\
\dim\ker\Phi(A)^* &=& \sum_{\ka_k>0}\Theta(\ro_k,\ka_k).\label{f5.ker2}
\eqn
Moreover, a pseudoinverse of $\Phi(A)$ is given by
\begin{equation}
\Psi(A\iv A_-)(I-\textstyle\frac{1}{2}\cH(DW))\cT(A_-\iv).
\end{equation}
\item[(b)]
Assume that an antisymmetric factorization of the
function $G(t)=\wt{A}\iv(t)WA(t)$ is given by
$G(t)=\wt{A}_+\iv(t)D(t)A_+(t)$, where $A_+\in G\cB_+\NN$ and
$D(t)$ is of the form (\ref{f3.D(t)}) with the
characteristic pairs (\ref{f3.pairs1}). Then
\bqn\label{f5.ker3}
\dim\ker\Psi(A)   &=& -\sum_{\ka_k<0}\Theta(-\ro_k,\ka_k),\\
\dim\ker\Psi(A)^* &=& \sum_{\ka_k>0}\Theta(-\ro_k,\ka_k).\label{f5.ker4}
\eqn
Moreover, a pseudoinverse of $\Psi(A)$ is given by
\begin{equation}
\cT(A_+\iv)(I-\textstyle\frac{1}{2}\cH(D\iv W))\Phi(A_+A\iv)
\end{equation}
\end{itemize}
\end{theorem}
\begin{proof}
The proof is based on Theorem \ref{t4.6} and Proposition \ref{p5.1}.
Because $\Phi(A)$ and $\Psi(A)$ are unitarily equivalent to
$\cM_W(A)$ and $\cN_W(A)$, respectively, the formulas for the
dimension of the kernel and cokernel follow immediately.
In order to show that the above expression are indeed pseudoinverses,
one can apply the C*-algebra isomorphism $\Xi$ introduced in
Proposition \ref{p5.1} to these operators. Using the formulas
stated there, one obtains
$$
\ba{rclcl}
\Xi &:& \Psi(A\iv A_-)(I-\textstyle\frac{1}{2}\cH(DW))\cT(A_-\iv)
&\mapsto&
\cN_W(A\iv A_-)(I-\textstyle\frac{1}{2}H(D)) T(A_-\iv),\\[.5ex]
\Xi &:& \cT(A_+\iv)(I-\textstyle\frac{1}{2}\cH(D\iv W))\Phi(A_+A\iv)
&\mapsto&
T(A_+\iv)(I-\textstyle\frac{1}{2}H(D\iv))\cM_W(A_+A\iv).
\ea
$$
The operators on the right hand side are exactly the
expressions (\ref{f4.pseu1}) and (\ref{f4.pseu2})
for the pseudoinverses of $\cM_W(A)$ and $\cN_W(A)$.
The observation that an operator $X\kr$ is a pseudoinverse
of an operator $X$ if and only if $\Xi(X\kr)$ is a pseudoinverse
of $\Xi(X)$ completes the proof.
\end{proof}

The corresponding invertibility criteria reads as follows
(compare Corollary \ref{c4.7}).

\begin{corollary}\label{c5.4}
Let $W$ be given by (\ref{f5.W}), and assume $A\in G\cB\NNN$. 
\begin{itemize}
\item[(a)] 
The operator $\Phi(A)$ is invertible if and only if the function
$F(t)=A(t)W\wt{A}\iv(t)$ admits an antisymmetric factorization with
characteristic pairs $(\ro_k,\ka_k)$
which are all contained in the set
\begin{equation}\label{f5.set1}
\Big\{\;(-1,-1),\,(-1,0),\,(1,0),\,(1,1)\;\Big\}.
\end{equation}
\item[(b)] 
The operator $\Psi(A)$ is invertible if and only if the function
$G(t)=\wt{A}\iv(t)WA(t)$ admits an antisymmetric factorization with
characteristic pairs $(\ro_k,\ka_k)$
which are all contained in the set
\begin{equation}\label{f5.set2}
\Big\{\;(1,-1),\,(-1,0),\,(1,0),\,(-1,1)\;\Big\}.
\end{equation}
\end{itemize}
\end{corollary}

\section{General Toeplitz + Hankel operators}

In this section we study Toeplitz + Hankel operators $T(a)+H(b)$
where no ``a priori'' relation between $a$ and $b$ is assumed.

It has been stated in Proposition \ref{p2.2} that the Fredholmness
of $T(a)+H(b)$ with $a,b\in L^\iy(\T)\NN$ implies $a\in G(L^\iy(\T)\NN)$.
Moreover, in the case where $a,b\in C(\T)\NN$ the necessary and sufficient
criteria for Fredholmness has been stated in Proposition \ref{p2.4}.

We are going to consider the case where $a,b\in\cB\NN$. It follows as before
from the inverse closedness of $\cB$ in $C(\T)$ that
$T(a)+H(b)$ with $a,b\in\cB\NN$ is Fredholm if and only if
$a\in G\cB\NN$. The dimension of the kernel and cokernel
in the case of Fredholmness reads as follows.

\begin{theorem}\label{t6.1}
Let $a\in G\cB\NN$, $b\in\cB\NN$, and $W$ be given by (\ref{f5.W}).
Introduce the functions
\bqn
A(t) &=& \left(\ba{cc} a(t)&b(t)\\ 0& I_N\ea\right)
\in\cB\NNN,\\
F(t) &=& A(t)W\wt{A}\iv(t)
\;\;=\;\;
\left(\ba{cc} b(t)\tilde{a}\iv(t) & a(t)-b(t)\tilde{a}\iv(t)\tilde{b}(t)\\
\tilde{a}\iv(t) & -\tilde{a}\iv(t)\tilde{b}(t)\ea\right)
\in\cB\NNN.\label{f6.Ft}
\eqn
If the characteristic pairs of the antisymmetric factorization
of $F(t)=A_-(t)D(t)\wt{A}_-\iv(t)$ are given by (\ref{f3.pairs1}), then
\bqn
\label{f6.ker1}
\dim\ker(T(a)+H(b))   &=& -\sum_{\ka_k<0}\Theta(\ro_k,\ka_k),\\
\label{f6.cok1}
\dim\ker(T(a)+H(b))^* &=& \sum_{\ka_k>0}\Theta(\ro_k,\ka_k).
\eqn
Moreover, if we write
\bqn
A\iv A_- \;\;=\;\;
\left(\ba{c} u_1\\ u_2\ea\right),\qquad
A_-\iv \;\;=\;\;
\Big(v_1,\;v_2\Big),
\eqn
with $u_1,u_2\in\cB^{N\times 2N}$ and $v_1,v_2\in\cB^{2N\times N}$, then
a pseudoinverse of $T(a)+H(b)$ is given by
\bqn\label{f6.pinv1}
\Big(T(u_1)+H(\tilde{u}_2)\Big)\Big(P-{\textstyle\frac{1}{2}}H(D)\Big)
T(v_1)
\eqn
\end{theorem}
\begin{proof}
The dimension the kernel and cokernel of the
Toeplitz + Hankel operator $T(a)+H(b)$, which is defined on
$(H^2)^N$, coincides with that of the operator
\bqn
X &=&
PM(a)P+PM(b)JP+Q,\nn
\eqn
which is defined on $(L^2)^N$. Now we write
\bqn
PM(a)P+PM(b)JP+Q &=&
\big(I-PM(a)Q-PM(b)JQ\big)
\big(PM(a)+PM(b)J+Q\big)\nn
\eqn
Because $I-Y=PM(a)Q+PM(b)JQ$ is nilpotent, the first expression on
the right hand side (i.e., the operator $Y$) is invertible. Hence we have
to determine the dimension of the kernel and cokernel of
\bqn
PM(a)+PM(b)J+Q,\nn
\eqn
which is just the singular integral operator $\Phi(A)$ with
$A(t)$ given as above. Now the result follows from
Theorem \ref{t5.4}.

As to the pseudoinverse, we first remark that a pseudoinverse of
$T(a)+H(b)$ is given by $PX\kr P$, where $X\kr$ is a pseudoinverse
of the above operator $X$. Since $X=Y\Phi(A)$, it follows that
$X\kr=(\Phi(A))\kr Y\iv$.
Hence $PX\kr P= P(\Phi(A))\kr P$ because
$P=Y\iv P$ as can easily be seen.
From Theorem \ref{t5.4}(a) we conclude that $(\Phi(A))\kr$ may be given by
$$
\Psi(A\iv A_-)(I-{\textstyle\frac{1}{2}}\cH(DW))\cT(A_-\iv)
$$
Using the definition of the operators occurring there, we obtain that this is
equal to
\bqn
&&\Big(I,J\Big)M(A\iv A_-)\left(\ba{c} P\\PJ\ea\right)
\left(I-\frac{1}{2}\Big(P,JP\Big)M(D)
\left(\ba{cc}QJ\\Q\ea\right)\right)
\Big(P,JP\Big)M(A_-\iv)\left(\ba{c}P\\PJ\ea\right)
\nn\\
&&\qquad=\quad
\Big(I,J\Big)M(A\iv A_-)\Big(P-\frac{1}{2}PM(D)JP\Big)
M(A_-\iv)\left(\ba{c}P\\PJ\ea\right).\nn
\eqn
Hence $PX\kr P=P(\Phi(A))\kr P$ equals
\bqn
\Big(P,PJ\Big)M\left(\ba{c} u_1\\u_2\ea\right)P
\Big(P-\frac{1}{2}H(D)\Big)
PM(v_1,v_2)\left(\ba{c}P\\0\ea\right),\nn
\eqn
which in turn is equal to the operator (\ref{f6.pinv1}).
\end{proof}

We want to emphasize that the matrices $u_k$ and $v_k$
are of size $N\times 2N$ and $2N\times N$, respectively, whereas $D$ is
of size $2N\times 2N$. The occurring Toeplitz and Hankel operators
are block operators of a corresponding size and their definition should
be obvious.

In the above theorem, we reduced the calculation of the dimensions and
pseudoinverse for $T(a)+H(b)$ to those of the singular integral operator
$\Phi(A)$. For reasons of symmetry, one should suspect that it can also be done
by reduction to the singular integral operator $\Psi(B)$.
We will establish the corresponding statement in the following theorem
for completeness sake. Of course, the corresponding result should
essentially be the same. How the assertions of
both of theorems are related with each other will be discussed afterwards.

\begin{theorem}\label{t6.2}
Let $a\in G\cB\NN$, $b\in\cB\NN$, and $W$ be given by (\ref{f5.W}).
Introduce the functions
\bqn
B(t) &=& \left(\ba{cc} a(t)&0\\ \wt{b}(t)& I_N\ea\right)
\in\cB\NNN,\\
G(t) &=& \wt{B}\iv(t)WB(t)
\;\;=\;\;
\left(\ba{cc} \tilde{a}\iv(t)\tilde{b}(t) &\tilde{a}\iv(t) \\
a(t)-b(t)\tilde{a}\iv(t)\tilde{b}(t) & -b(t)\tilde{a}\iv(t)\ea\right)
\in\cB\NNN.\label{f6.Gt}
\eqn
If the characteristic pairs of the antisymmetric factorization
of $G(t)=\wt{B}\iv_+D(t)B_+(t)$ are given by (\ref{f3.pairs1}), then
\bqn
\label{f6.ker2}
\dim\ker(T(a)+H(b))   &=& -\sum_{\ka_k<0}\Theta(-\ro_k,\ka_k),\\
\label{f6.cok2}
\dim\ker(T(a)+H(b))^* &=& \sum_{\ka_k>0}\Theta(-\ro_k,\ka_k).
\eqn
Moreover, if we write
\bqn
B_+\iv \;\;=\;\;
\left(\ba{c} w_1\\ w_2\ea\right),\qquad
B_+B\iv \;\;=\;\;
\Big(y_1,\;y_2\Big),
\eqn
with $w_1,w_2\in\cB^{N\times 2N}$ and $y_1,y_2\in\cB^{2N\times N}$, then
a pseudoinverse of $T(a)+H(b)$ is given by
\bqn\label{f6.pinv2}
T(w_1)\Big(P-{\textstyle\frac{1}{2}}H(D\iv)\Big)
\Big(T(y_1)+H(y_2)\Big)
\eqn
\end{theorem}
\begin{proof}
As before, the dimension of the kernel and cokernel of
$T(a)+H(b)$ coincides with that of the operator
$X=PM(a)P+PM(b)JP+Q$. Now we write
\bqn
PM(a)P+PM(b)JP+Q &=&
\big(M(a)P+M(b)JP+Q\big)
\big(I-QM(a)P-QM(b)JP\big)\nn
\eqn
Because $Y'=QM(a)P+QM(b)JP$ is nilpotent, the last expression on
the right hand side is an invertible operator. Hence we are led
to the dimension of the kernel and cokernel of
\bqn
M(a)P+M(b)JP+Q,\nn
\eqn
which coincides with the singular integral operator $\Psi(B)$ where
$B(t)$ is given as above. Now the result follows from
Theorem \ref{t5.4}.

Again, a pseudoinverse of
$T(a)+H(b)$ is given by $PX\kr P$. Since $X=\Psi(B)$Y', it follows that
$X\kr=(Y')\iv(\Psi(B))$.
Hence $PX\kr P= P(\Psi(B))\kr P$ because
$P=P(Y')\iv$ as can easily be seen.
From Theorem \ref{t5.4}(b) we conclude that $(\Psi(B))\kr$ may be given by
$$
\cT(B_+\iv)(I-{\textstyle\frac{1}{2}}\cH(D\iv W))
\Phi(B_+ B\iv)
$$
We obtain that this is equal to
\bqn
\hspace{-3ex}
&&\Big(P,JP\Big)M(B_+\iv)\left(\ba{c} P\\PJ\ea\right)
\left(I-\frac{1}{2}\Big(P,JP\Big)M(D\iv)
\left(\ba{cc}QJ\\Q\ea\right)\right)
\Big(P,JP\Big)M(B_+B\iv)\left(\ba{c}I\\J\ea\right)
\nn\\
\hspace{-3ex}
&&\qquad=\quad
\Big(P,JP\Big)M(B_+\iv)\Big(P-\frac{1}{2}PM(D\iv)JP\Big)
M(B_+B\iv)\left(\ba{c}I\\J\ea\right).\nn
\eqn
Hence $PX\kr P=P(\Psi(B))\kr P$ equals
\bqn
\Big(P,0\Big)M\left(\ba{c} w_1\\w_2\ea\right)P
\Big(P-\frac{1}{2}H(D\iv)\Big)
PM(y_1,y_2)\left(\ba{c}P\\JP\ea\right),\nn
\eqn
which in turn is equal to the operator (\ref{f6.pinv2}).
\end{proof}

Now we discuss the question of how the statements of the preceding
two theorems are related with each other. At first glance, formulas
(\ref{f6.ker1}) and (\ref{f6.cok1}) seem to contradict
(\ref{f6.ker2}) and (\ref{f6.cok2}) but this is just because the same notation
has been used for different factors $D(t)$ with different characteristic pairs.
In these theorems we start with the antisymmetric factorizations
of certain functions $F(t)$ and $G(t)$:
$$
F(t)\;\;=\;\;A_-(t)D^{(1)}(t)\wt{A}_-\iv(t),
\qquad
G(t)\;\;=\;\;\wt{B}_+\iv(t)D^{(2)}(t)B_+(t).
$$
Assume that the notation of the characteristic pairs is given by
$D^{(1)}(t)=\diag(\ro_i^{(1)}t^{\ka_i^{(1)}})$ and
$D^{(2)}(t)=\diag(\ro_i^{(2)}t^{\ka_i^{(2)}})$.
The functions $F(t)$ and $G(t)$ are given by (\ref{f6.Ft}) and
(\ref{f6.Gt}), from which it follows that
$$
G(t)\;\;=\;\;
- \left(\ba{cc}0&1\\-1&0\ea\right) F(t)
\left(\ba{cc}0&-1\\1&0\ea\right)
$$
Note that the constant matrices on the right hand side are the inverses of
one another. The last identity is the reason that from an
antisymmetric factorization of $F(t)$
one can immediately obtain an antisymmetric factorization of $G(t)$ and vice
versa. More precisely, if we are given an antisymmetric factorization of $F(t)$,
then an antisymmetric factorization of $G(t)$ is given with the factors
$$
B_+(t)\;\;=\;\;\wt{A}_-\iv(t)\left(\ba{cc}0&-1\\1&0\ea\right)
\quad\mbox{ and }\quad
D^{(2)}(t)\;\;=\;\;-D^{(1)}(t).
$$
This also shows that
the construction of an antisymmetric factorization for $F(t)$ and
$G(t)$ is essentially the same problem. Moreover, if the characteristic
pairs are ordered ``appropriately'', we may conclude that
$(\ro_k^{(2)},\ka_k^{(2)})=(-\ro_k^{(1)},\ka_k^{(1)})$ for all $k$.
This implies that formulas (\ref{f6.ker1}) and (\ref{f6.cok1}) 
indeed coincide with (\ref{f6.ker2}) and (\ref{f6.cok2}).

It is, however, not clear at this point whether formulas
(\ref{f6.pinv1}) and (\ref{f6.pinv2}) for the pseudoinverse are the same.
Observe that pseudoinverses are in general not unique. Of course, if the
pseudoinverses are inverses, then they automatically have to be the same.

\section{More general singular integral operators}

In this section we consider a more general class of
singular integral operators. The operators $\Phi(A)$ and $\Psi(A)$
just represent special cases of this class of operators
(see (\ref{f5.124}) and (\ref{f5.125})). These more general singular
integral operators are operators of the form
\bqn\label{f7.SIOF}
&&PM(a_1)P+PM(b_1)JQ+QM(\tilde{c}_1)JP+QM(\tilde{d}_1)Q+\\[1ex]
&&PM(a_2)JP+PM(b_2)Q+QM(\tilde{c}_2)P+QM(\tilde{d}_2)JQ,
\eqn
where $a_i,b_i,c_i,d_i\in L^\iy(\T)\NN$, $i=1,2$. Introducing the
functions $A,B\in L^\iy(\T)\NNN$ and the constant $W\in\C\NNN$ by
\bqn
A=\left(\ba{cc} a_1&b_1\\c_1&d_1\ea\right),\quad
B=\left(\ba{cc} a_2&b_2\\c_2&d_2\ea\right),\quad
W=\left(\ba{cc} 0&I_N\\I_N&0\ea\right),
\eqn
it is easily seen by help of formulas (\ref{f5.cT})
and (\ref{f5.cH}) that the operator (\ref{f7.SIOF})
equals
\bqn\label{f6.SIOF2}
\cT(A)+\cH(BW).
\eqn

The following result is an immediate consequence of Proposition \ref{p2.2}
and Proposition \ref{p2.4} in connection with Proposition \ref{p5.1}.
It is also a generalization of Proposition \ref{p5.2}.

\begin{proposition}\label{p7.1}
Let $A,B\in L^\iy(\T)\NNN$.
\begin{itemize}
\item[(a)]
If $\cT(A)+\cH(BW)$ is Fredholm, then $A\in G(L^\iy(\T)\NNN)$.
\end{itemize}
Let $A,B\in C(\T)\NNN$.
\begin{itemize}
\item[(b)]
$\cT(A)+\cH(BW)$ is Fredholm if and only if $A\in G(C(\T)\NNN)$.
\end{itemize}
Moreover, if this is true, then $\ind(\cT(A)+\cH(BW))=-\wind\det A$.
\end{proposition}

In fact, the mapping $\Xi$ defined in Proposition \ref{p5.1} sends the
operator $\cT(A)+\cH(BW)$ into the Toeplitz + Hankel operator
$T(A)+H(B)$. In the case where $A\in G(\cB\NNN)$ and $B\in\cB\NNN$,
these operators are Fredholm, and formulas for the dimension of the
kernel and cokernel can be obtained by help of the results of the previous
section.

\begin{theorem}\label{t7.2}
Let $A\in G(\cB\NNN)$ and $B\in\cB\NNN$. Introduce the function
$F$ by
\bqn
F &=& 
\left(\ba{cc} B\wt{A}\iv & A-B\wt{A}\iv\wt{B} \\
\wt{A}\iv & -\wt{A}\iv\wt{B}\ea\right)\in G(\cB^{4N\times 4N}).
\eqn
Then $F$ admits an antisymmetric factorization. If the characteristic
pairs are denoted by $(\ro_k,\ka_k)$, $k=1\dots 4N$, then
\bqn
\dim\ker(\cT(A)+\cH(BW)) &=& -\sum_{\ka_k<0}\Theta(\ro_k,\ka_k),\\
\dim\ker(\cT(A)+\cH(BW))^* &=&\sum_{\ka_k>0}\Theta(\ro_k,\ka_k).
\eqn
\end{theorem}

It is also possible to establish formulas for the pseudoinverses of
$\cT(A)+\cH(BW)$. We leave these details to the reader.

\end{document}